\documentclass{article}
\usepackage{graphicx,verbatim}
\usepackage{psfrag}
\usepackage{amsfonts}
\usepackage{amssymb}
\usepackage{amsmath}
\usepackage{lscape}
\usepackage{epsf}
\usepackage{bm}
\newtheorem{thm}{Theorem}

\newtheorem{lemma}[thm]{Lemma}

\newcommand\C{\mathbb{C}}
\newcommand\bN{\mathbb{N}}

\newcommand\R{\mathbb{R}}
\def\rho{\varrho}
\newcommand\bI{\mathbb{I}}
\newcommand\bR{\R}

\newcommand\N{\mathbb{N}}

\newcommand\hmc{{\bf {C}}}
\newcommand\ma{{\bf {R}}}
\newcommand\gmp{{{\bf {P}}}}

\title{Width of homoclinic zone for quadratic maps.}

\author{V. Gelfreich $\dag$, V. Naudot $^\star$
\\
$ \ $
\\
Mathematics Institute\\
University of Warwick\\
Coventry CV4 7AL, UK
\\
$ \ $
\\
email: $\dag$ V.Gelfreich@warwick.ac.uk, 
\\
$ \ \ \ \ \ \  $  
\footnote{Corresponding author
} $ \ $V.Naudot@warwick.ac.uk}

\begin{document}
\maketitle
\vfill
\eject
\par\vskip10pt
\begin{abstract}
We study several families of planar quadratic diffeomorphisms 
near a Bogdanov-Takens bifurcation.
For each family, the associate bifurcation diagram 
can be deduced from the interpolating flow.
However, a zone of chaos confined between two  lines of homoclinic bifurcation
that are exponentially close to one-another 
 is
observed.
The goal of this  paper is to test numerically an accurate asymptotic 
expansion
for 
 the width of this
chaotic zone for different families.

\end{abstract}

\par\vskip10pt
Mathematics Subject Classification:
37D45, 37E30, 37G10.

\section{Introduction}
In this paper we study 
homoclinic bifurcations in the 
unfolding of a diffeomorphism
near a fixed point of Bogdanov-Takens type. 
To begin with, 
we consider a planar diffeomorphism 
$F :{\bR}^2 \rightarrow {\bR}^2$
with the origin as a fixed point and where
\begin{eqnarray*}
dF(0,0) & = & {{\bI}}{\mathrm {d}} +{\mathrm{ N}}
\end{eqnarray*}
where ${\mathrm {N}} {\not \equiv} 0$ is nilpotent.
The origin is said 
to be a fixed point of {\it Bogdanov-Takens}
 type.
This latter terminology is more known for a
singularity of a vector field $X$
with linear part
having double zero 
eigenvalues and a non vanishing nilpotent part.
Since this singularity is of codimension 2,. i.e., is twice degenerate, 
a 
generic unfolding 
will depend on two parameters say $(\mu,\nu)$.
In the case of a vector field,
such unfolding has been studied in \cite{a,Takens1974} and for maps
in \cite{BroerRS1996, brs2}.
For completeness, the 
corresponding bifurcation diagram is 
 revisited in Figure 1 on the left:      
a curve
of homoclinic bifurcation emanates from the  
origin,
 below a curve of 
Hopf bifurcation, see \cite{dtbt} for 
the terminology and more details.
For parameters located between these two curves,
 the corresponding dynamics possesses a stable limit cycle.
Finally, for parameter on the ordinate $\{ \mu=0 \}$, a saddle node
occurs, 
see also
\cite{dtbt} for more details.

The Bogdanov-Takens bifurcation plays an important role in dynamical systems, 
for instance
from the bifurcation theoretical point of view. Given any dynamical systems depending on a parameter,
 the structure of the bifurcation set can be often understood by the presence of 
several high codimension points which act as organising centres. Knowing the presence of
(degenerate or not degenerate)  Bogdanov Takens points initiate the searches for 
subordinate bifurcation sets such as Hopf bifurcation sets or homoclinic bifurcation sets.
In this paper, we consider a nondegenerate Bogdanov Takens point.

For the map $F$, an unfolding theory is developed 
in
\cite{BroerRS1996,brs2}.  It is very similar to the case of a flow.
To be more precise, 
any unfolding 
\begin{eqnarray*}
{ F}_{\mu,\nu} & : & 
{\bR}^2 \rightarrow {\R}^2, \ (x,y) \mapsto (x_1,y_1)  
\end{eqnarray*}
of the map $F$
(where $(\mu,\nu) \in {\bR}^2$ and $F=F_{0,0}$)
can be  
 embedded  into a nonautonomous and periodic
family of vector fields $X_{\mu,\nu}$. 
The diffeomorphism coincides with the time 1 map of 
that vector field, see also \cite{Takens1974}.
Using an averaging theorem \cite{nei}
 the dependence on 
time is removed 
to exponentially 
small terms.
Moreover, one can show that $F_{\mu,\nu}$ 
 is formally interpolated by an autonomous vector field
${\tilde X}_{{ \mu}, { \nu}}$, see \cite{Gelfreich2003}.
This latter can be used to
study the bifurcations of fixed points of $F_{\mu,\nu}$.
Both approaches move the difference between these two types of bifurcations 
beyond all
algebraic order.

 Although all the Taylor coefficients of
$ {\tilde X}_{{ \mu}, {\nu}}$  can be written, 
there is no reason to expect convergence 
of the corresponding series, 
since 
the dynamics 
for a planar 
diffeomorphism can be much richer 
than the dynamics of a planar vector field.
In the real analytic theory, this difference is 
exponentially small 
\cite{BroerR2001, Gelfreich2003}.

\begin{figure}[ht33333]
\includegraphics[scale=0.60]{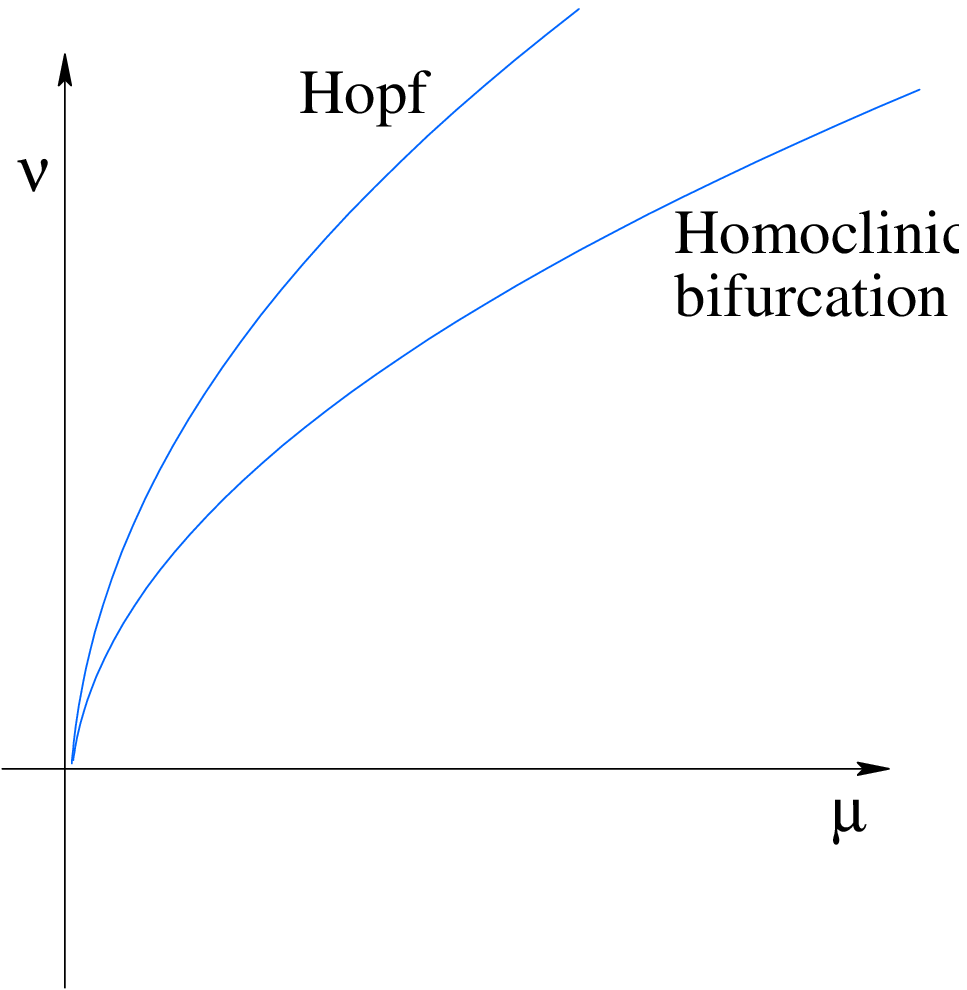}
\includegraphics[scale=0.60]{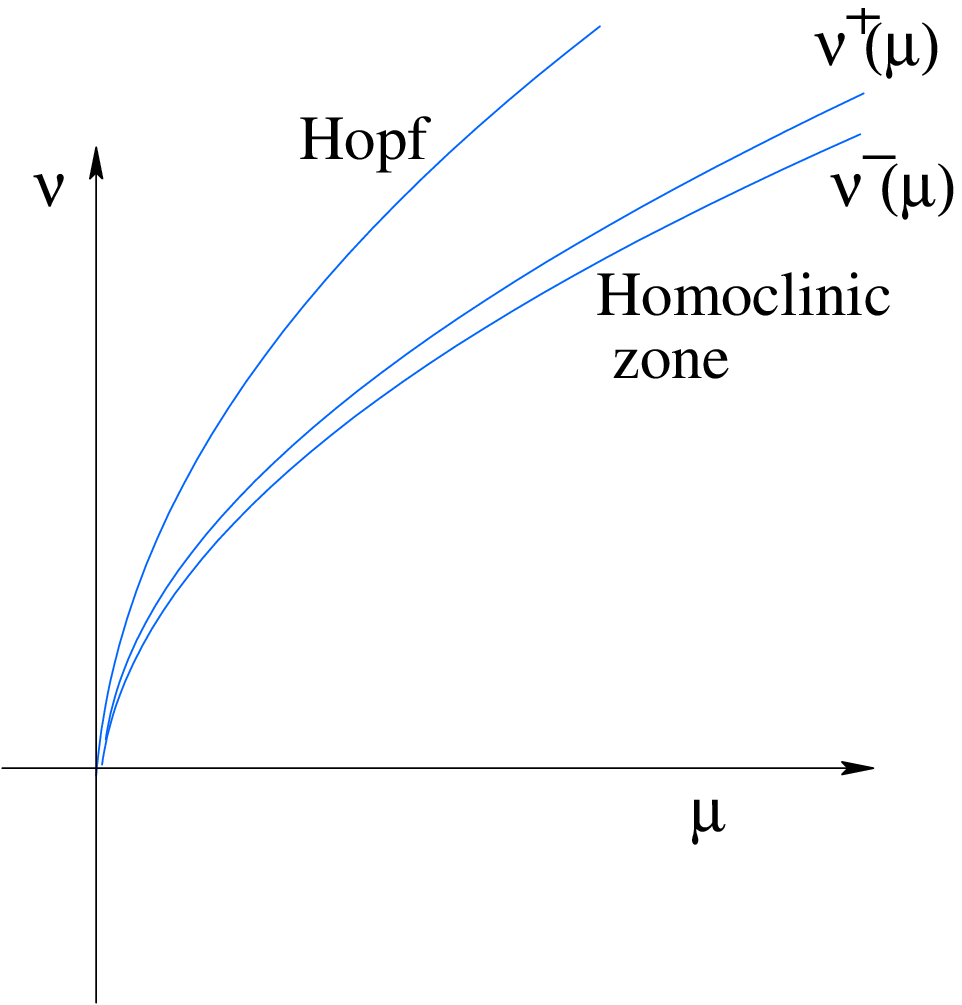}
\caption{The Bogdanov Takens bifurcation for  a flow (left) 
and for a diffeomorphism
(right).}
\end{figure}\label{fig:BT}

As we said above,
for diffeomorphisms, the bifurcation diagram 
(figure 1, on the right)
is essentially the same. However,
there is no reason to expect a 
single homoclinic curve, 
since a homoclinic orbit 
may be transverse  
and therefore persists.
We observe a separatrices splitting and
instead of a single homoclinic curve, 
one observes 
two curves $\nu^+(\mu) $ and $\nu^-(\mu)$ 
respectively  
corresponding to the first and  the last homoclinic tangency.
If a parameter $(\mu, \nu)$ is (strictly) located in region 
between those two curves,
then the map $F_{\mu,\nu}$ possesses transverse homoclinic trajectories. 
On the lower and upper boundary
the homoclinic connexion becomes non-transverse.
Understanding the width of this region is the main goal of this paper.

Before going any further, we set the following preliminaries.
Without lost of generality and
up to an analytic change of coordinates one has:
\begin{eqnarray}\label{unf}
 x_1 = x+y,  \   y_1 = y + f_{\mu,\nu}(x,y), 
\end{eqnarray}
where
$$
f_{0,0}(0)
=0=
\frac{\partial f_{0,0}}{\partial x}(0,0)
=\frac{\partial f_{0,0}}{\partial y}(0,0).
$$
We shall assume that
\begin{eqnarray}
\label{q2}
\displaystyle\frac{\partial^2 f_{0,0}}{\partial x^2}(0,0) {\not =} 0.
\end{eqnarray}
By the implicit function theorem, there exists
${\tilde x}({\mu,\nu}) $ such that
$$
\frac{\partial f_{\mu,\nu}}{\partial x}({\tilde x}({\mu,\nu}),0) \equiv 0
.$$
Applying a 
conjugacy of the form $x ={\bar x}+ {\tilde x}_{\mu,\nu}, \ y ={\bar y}$
(and after removing the bars)
amounts to writing
\begin{eqnarray}\nonumber
f_{\mu,\nu} (x,y)
&  =&   
-b_{00}(\mu,\nu) + b_{20}(\mu,\nu)x^2   + b_{01}(\mu,\nu) y
\\ \label{unfold}
 & + & b_{11}(\mu,\nu)  x y+
{\mathrm {h.o.t}}(x,y)
\end{eqnarray}
${\mathrm {h.o.t}}(x,y)$ stands for the higher order terms in $x$ and $y$.
From
(\ref{q2}) $b_{20}(0,0) {\not =} 0$.
By a linear rescaling in the variables $(x,y)$,
we can fix $b_{20}(\mu,\nu) \equiv 1$. Furthermore, we 
 put $b_{11}(0,0)=\gamma$
and assume that the map
 $$ 
(\mu,\nu) \mapsto (-f_{\mu,\nu}(0,0),
\displaystyle\frac{\partial f_{\mu,\nu}}{\partial y}(0,0))=
(-b_{00}(\mu,\nu), b_{01}(\mu,\nu))
$$
is a local diffeomorphism near $(0,0)$. From now on, we shall consider
$(b_{00}, b_{01})$ as our parameters and rename them (again) by
$(\mu,\nu)$ i.e., write
$(b_{00}, b_{01}) = (\mu,\nu)$.

\vskip10pt
\noindent
In \cite{BroerRS1996} it was shown that 
\begin{eqnarray}\label{homo-line}
\nu^\pm (\mu) & = & \frac{5}{7}(\gamma-2){\sqrt{\mu}} + 
{\cal O}(\mu^\frac{3}{4}).
\end{eqnarray}
In \cite{Gelfreich2003}
 the following formula is proposed:
\begin{eqnarray}\label{gelf}
\nu^+ (\mu)-\nu^-(\mu) & = & { {\mathbf{\Theta}}}_\gamma
K(\mu,\gamma-2) +{\cal O}(\mu^{1/4}\log\mu),
\end{eqnarray}
where
\begin{eqnarray}\label{kmugamma}
K(\mu,{\tilde \gamma}) =
\frac{5}{6\sqrt{2}\mu^{\frac{5}{4}}}\cdot e^{-\sqrt{2}\pi^2/\sqrt[4]{\mu}}
e^{-6\pi^2{\tilde \gamma}/7}
\end{eqnarray}
is referred to as the `leading part' of the width
and  
${{\mathbf{\Theta}}}_\gamma$ is 
an analytic invariant of 
the map $F_{0,0}$ called a `splitting constant', see \cite{GN}. 

\vskip10pt
The goal of this paper is to establish, numerically, a more 
accurate formula for the width of the homoclinic zone 
$\nu^+ (\mu)-\nu^-(\mu)$.
The existence of asymptotic 
expansions for the width of the homoclinic zone
is  unknown so far.
 Furthermore,
if it does
exist, it is very hard to compute analytically. The difficulty here
comes from the fact that the normal form of the map
coincides with that of the  time 1-map of a vector field.
Therefore the difference between  the flow and the map
is pushed beyond any algebraic order. 
In the nearly integrable context,  a 
polynomial asymptotic expansion 
for the splitting of the separatrices
is proposed in \cite{ros1}: 
the author considers the perturbation of a Hamiltonian (elliptic) billiard.
The system depends on  a perturbation parameter $\varepsilon \geq 0$,  
a hyperbolicity parameter
$h > 0$ 
and admits
 four separatrices, which break up when $h>0$.
In this special case, 
the author proposed an asymptotic expansion 
for the area
of the main lobes of the resulting turnstile
that 
takes the form of a power series (with even terms) 
in
$\varepsilon$.
See \cite{GS, ros2, lev-tab} 
for more references on the computation of 
separatrices splitting.

In this paper our approach is somehow experimental.
We study examples and present  
strong numerical evidence for 
the following expansion of the width of the homoclinic zone
\begin{eqnarray}
\label{teo}
\nu^+(\mu)-\nu^-(\mu) 
 \asymp  
K(\mu,\gamma-2) 
\sum_{ k\geq 0}\sum_{0\leq j\leq [\frac{k}{2}]}{\tilde c}_{k,j}\mu^{k/4}\log^{j} \mu
\end{eqnarray}
where $\displaystyle[\frac{k}{2}]$ stands for the integer part of $k/2$ and 
$K$ is given by (\ref{kmugamma}).

\par\vskip10pt
\noindent

{\bf Remarks}:\par\vskip10pt

(a) Observe that (\ref{teo}) 
 is a double series with logarithmic terms and numerically,
for such an expansion, we do not know any
efficient techniques to compute  the corresponding coefficients with a large
precision.
However, 
our numerical experiments showed that 
$\log(\nu^+(\mu)-\nu^-(\mu))$ has a simpler asymptotics 
expansion.
than $\nu^+(\mu)-\nu^-(\mu)$ itself. 
 More precisely we have
\begin{eqnarray}
\label{teo-log}
\log\biggl(\nu^+(\mu)-\nu^-(\mu)\biggr) 
 \asymp   \displaystyle
\log\biggl(K(\mu,\gamma-2)\biggr) + 
\sum_{k\geq 1} {m}_k \mu^{\frac{k}{4}} +\log \mu \sum_{k\geq 1} 
n_{k} \mu^{{\frac{k}{2}}}. 
\end{eqnarray}
One easily checks that formula (\ref{teo}) follows from (\ref{teo-log}) and that
the ${\tilde c}_{k,j}$'s depend on the $m_k$'s and the $n_k$'s.
Note that the asymptotic series (\ref{teo-log})
does not involve a double summation and therefore the corresponding
coefficients can be computed with a much higher precision.

\par\vskip10pt
(b) Logarithmic terms may vanish, this occurs for instance in the case
of the H\'enon map, see next section for more details.

\par\vskip10pt
(c) 
From the numerical data, we are able to guess
a simple analytic expression for 
the first logarithmic term in (\ref{teo-log}). More precisely
we have
$$
n_1 = -\biggl(\frac{6(\gamma-2)}{7\sqrt{2}}\biggr)^2,
$$
which is valid for all families studied in this article.
The paper is organised as follows. 
We shall consider three
different families that satisfy, 
(up to appropriate 
smooth
changes of coordinates)
 the setting above with different non linear terms.
As a result of our experiments, for each family we shall 
state
the asymptotics for the width of the homoclinic zone,
confirming formula (\ref{teo-log}).
Looking for the width of the homoclinic zone amounts to
fixing one parameter, say $\mu$ in the unfolding (\ref{unfold}),
and find the values of the second parameter, say $\nu$, 
 for which the system
admits  a first and a last homoclinic tangency.
We say that $\mu$ is the main `parameter'  and $\nu$ is 
the `slave' parameter.
In section 3 we briefly present the strategy to follow. The rest of the section
is devoted to 
 the computation
of the invariant (stable and unstable) 
manifolds at the saddle point. The {\it splitting function} 
which is a  key ingredient of the techniques 
is presented.
Indeed,
primary homoclinic orbits are  
in one to one correspondence with zeroes of the splitting function. 
Therefore, the first and the last homoclinic  tangencies
will correspond to double zeroes of the splitting functions.
Moreover, the splitting function is periodic, with exponentially decreasing harmonics
and is well approximated by the {\it splitting determinant}.  
With a good precision, 
computing the width of the zone amounts to the computation of the
first two harmonics of
the splitting function 
and their dependence with respect to the slave parameter 
(the 
main parameter being fixed).
For each family, we compute the width of the homoclinic zone for several
hundreds values of the main parameter $\mu$
and collect the results
in a set of renormalised data. 
In the next step, 
 the coefficients in (\ref{teo-log}) (considered as an ansatz)
are extracted 
by  interpolation
techniques.
The remaining part of the paper 
is devoted to the verification of  
 the validity of 
our results.
More precisely, we test the ansatz (\ref{teo-log})
and we find how precise our data for the width of the homoclinic zone
should be in order to produce reliable results for the coefficients
of the asymptotic expansion.
Finally, 
the constant coefficient of the expansion should coincide with
the splitting constant \cite{Gelfreich2003}: 
following the procedure 
developed in \cite{GN}, we compare these constants with
the constant coefficients of the expansions.

\section{Main results} 

Before presenting our main results, we
first introduce the following notions.

\subsection{Asymptotic sequences and expansions}
 \par\vskip10pt
\noindent
Let $\varepsilon_0>0$ be given and
let 
\begin{eqnarray*}
{\tilde {\cal S}} =\{f_0,f_1,\ldots,f_n,\ldots \}
\end{eqnarray*}
where $f_0 \equiv 1$ and for each integer $i>0$,
$f_i:(0,\varepsilon_0)\rightarrow {\bR} $ is a smooth positive 
function 
such that
 $$\displaystyle{\lim_{x\to 0^+} \frac{f_{i+1}(x)}{f_i(x)}}=0,$$
or in other words 
$f_{i+1}(x)=o(f_{i}(x))$.
 Such a family ${\tilde {\cal S}}$ is called an {\it asymptotic sequence}.
In this paper we shall consider
the 
following  asymptotic sequences
\begin{eqnarray}\label{polyn}
{\tilde {{\cal P}}} & = & \{1,x,x^2\ldots,x^n,\ldots \}
\end{eqnarray}
that is $f_i(x) =x^i$
 and
the Dulac asymptotic sequence \cite{mar2}:
$$
{\tilde {\cal D}}= 
\{1,x,x^2\log x,
x^2,x^3,x^4\log x,x^4,\ldots, x^{2n}\log(x),x^{2n},x^{2n+1},\ldots \} 
$$
that is
for all integer $n \geq 0$
\begin{eqnarray}\label{dul}
f_{3n}(x)=x^{2n}, \ f_{3n+1}(x)=x^{2n+1}, \ f_{3n+2}(x)=x^{2n+2}\log(x).
\end{eqnarray}
\par\vskip10pt
\noindent

Let ${\phi }: (0,\varepsilon_0) 
\rightarrow {\bR}$ be a smooth function. We say that
\begin{eqnarray}\label{def-b}
{\phi }(x) & \asymp & \sum_{n \in {\bN}} \alpha_n f_n(x)
\end{eqnarray}
is an asymptotic expansion of ${\phi }$ at $0$
(where the $\{f_n\}_{n \in {\bN}}$ is an asymptotic sequence and 
all $\alpha_n$'s are real)
if for all integer $n$,

\begin{eqnarray*}
{\phi}(x) -  {\phi }^{\{ n\}} (x) = {\cal O}(f_{n+1}(x)), \ {\mathrm{where}} 
\ {\phi }^{\{ n\}} (x) = \sum_{i=0}^n \alpha_i f_i(x).
 \end{eqnarray*}
When looking at expansions 
of the form (\ref{def-b})
no convergence is implied 
and often the $a_i$'s are Gevrey-1, i.e.,
\begin{eqnarray}\label{gevrey}
\exists\   M >0, \ r>0, \ {\mathrm {such}} \ {\mathrm {that }} 
\ \forall k\geq 0, \  |\alpha_k| \leq  M k!/r^k.
\end{eqnarray}

\par\vskip10pt


 \subsection{Quadratic family}
Our first example is the Quadratic map
\begin{eqnarray}  
\label{rfam}
{\bf Q}={\bf Q}_{\mu, \nu,\gamma} & :& {\bR}^2 \rightarrow {\bR}^2 \\   \nonumber
(x,y) & \mapsto  & (x+y, y+x^2-\mu+ \gamma x y +\nu y) 
\end{eqnarray}
Observe that ${\bf Q}$  mimics 
the unfolding (\ref{unf}) i.e., 
that takes the form of
 ${F}_{\mu,\nu}$ and ignores the higher order terms.
We normalise 
the width of the 
homoclinic zone associated to the Quadratic family
by defining
\begin{eqnarray*}
S_\gamma(\mu)
& = &\frac{\nu^+(\mu)-\nu^-(\mu)}{K(\mu, \gamma-2)}.  
\end{eqnarray*}
Within precision 
of our computations we observe 
\begin{eqnarray}
\label{s-r-diverge}
\log S_\gamma (\mu) & \asymp &
\sum_{k\geq 0}M_k(\gamma) \mu^{k/4} +
\log \mu \sum_{k\geq 1} N_{k}(\gamma) \mu^{k/2}, 
\end{eqnarray}
where
$M_k(\gamma)$ and $N_k(\gamma)$ 
are real coefficients which depend on the parameter
$\gamma$.
Comparing with (\ref{gelf}), we see that
\begin{eqnarray*}
\exp(M_0(\gamma)) \equiv {{\mathbf{\Theta}}}_\gamma 
\end{eqnarray*}
is the
splitting constant
associated
with ${\bf Q}_{0,0,\gamma}$. Moreover, as 
we announced in the previous section, we
have
$$
N_1(\gamma) \equiv -\biggl(\frac{6(\gamma-2)}{7\sqrt{2}}\biggr)^2.$$
For each value of $\gamma$, 
the $M_k$'s and $N_k$'s 
can be computed with a very high precision, see Table 1 
for illustration.
Formula (\ref{s-r-diverge}) is verified for the 
76 first coeficients: $M_k$,  $k=0,\ldots,50$ and 
$N_\ell$, $\ell=1,\ldots,25$.
Although the precision decreases almost linearly as $k$ and $\ell$ increase, 
the 76 
first coefficients can be computed with 60 correct digits.
To compute these first coeficients, we need to compute the 
width of the homoclinic zone with at least 200 correct digits, 
see section 4.2 for more details.

Even if we can propose an analytic expression for $N_1(\gamma)$,
we have not been able to guess analytic expressions for the other coefficients
$N_k$ and $M_k$.

 \subsection{Bogdanov family}
Our second example is the
Bogdanov map \cite{arrow, ar2, Bogdanov1975}.
\begin{eqnarray*}  
{\bf B}={\bf B}_{a, b,{\tilde \gamma}} & :& {\bR}^2 \rightarrow {\bR}^2 \\  
   \nonumber
(x,y) & \mapsto  & 
(x+y+x^2 +
{\tilde \gamma} x y +a x +b y ,y+x^2 +{\tilde \gamma} x y +a x +b y). 
\end{eqnarray*}
The Bogdanov map, see for example \cite{a2,b2}, 
is the Euler map of a two-dimensional system of ordinary 
differential equations. 
In \cite{arrow} 
 Arrowsmith 
studied
the bifurcations and basins of attraction  and 
showed the existence  of mode locking, Arnold tongues, 
and chaos, see also \cite{ar2} for more details.

For this map the saddle point is located at the origin. 
This map can be transformed to the form (\ref{unfold}).
Indeed, let
$$
u= x - a/2, \ v= y+(x - a/2)^2 +{\tilde \gamma} (x-a/2) y  +a (x-a/2) +b y.
$$
We retrieve the map (2) and higher order terms (\ref{unfold}) 
by putting
$$
\nu=a+b -({\tilde \gamma}+2)\frac{a}{2}, \ \gamma =  {\tilde \gamma}+2, \ \mu=a^2/4,
$$
and 
$$
f_{\mu,\nu}= (x+y)^2-\mu +\gamma y^2. 
$$
The parameter $a$ is chosen to be the 
main parameter and $b$ the slave parameter.
From (\ref{homo-line}), the Bogdanov map admits a homoclinic zone near 
the line
$$
b^{\pm}(a) = \frac{6}{7}a {\tilde \gamma} +{\cal O}(a^{3/2}).
$$
The normalised width takes the form
\begin{eqnarray*}
 {\tilde S}_{\tilde \gamma} (a)
& = &  
\frac{b^+(a)-b^-(a)}{K(a^2/4, {\tilde \gamma})}.
\end{eqnarray*}
Similarly to the Quadratic family, 
our experiments showed that 
$\log {\tilde S}_{\tilde \gamma} $ satisfies
the following 
asymptotics:
\begin{eqnarray}
\label{bog-diverge-henon}
\log {\tilde S}_{\tilde \gamma} (a) & \asymp & 
\sum_{k\geq 0}A_k({\tilde \gamma}) a^{k/2} +\log a
\sum_{k\geq 1} B_{k}({\tilde \gamma}) 
a^{k},  
\end{eqnarray}
where
$A_k({\tilde \gamma})$ 
and 
$B_k({\tilde \gamma})$ are real coefficients
which depend on the parameter
$\gamma$. Comparing with (\ref{gelf}), we see that
 $$\exp(A_0({\tilde \gamma})) \equiv { {\mathbf{\Theta}}}_{\tilde \gamma}$$
is the
splitting constant
associated
with ${\bf B}_{0,0,{\tilde \gamma}}$.
Moreover, we observe numerically that 
$B_1({\tilde \gamma}) \equiv -(6{\tilde \gamma}/{7})^2$.
\vskip10pt

In Table 1, we provide typical results  for our computation
for the Quadratic and Bogdanov maps.
Although the first 20 coefficients do not show a tendency to grow
rapidly, we conjecture the series (\ref{s-r-diverge}) and
(\ref{bog-diverge-henon}) diverge and belong to the Gevrey-1 class (\ref{gevrey}), 
compare with \cite{GS}.
\par\vskip10pt
\begin{table}\label{T1}
\begin{tabular}{|c|c|c|}
\hline  
 coef. & scale & value  
\\ \hline
$ A_0$ &  $ 1 $ & ${\tt 61.26721889}$     \\ \hline
$A_1$ &  $a^{1/2} $  & ${\tt -29.82701974 }$    \\ \hline
$B_1$ & $ a\log a $  & ${\tt -6.612244898 } $   \\ \hline
$A_2$ & $ a $ & ${\tt 5.824479250  }$ \\ \hline
$A_3$ & $a^{3/2}$ & ${\tt 17.41183781  }$\\ \hline
$B_2$ & $a^2\log a $ &${\tt  5.649967276 }$  \\ \hline
$A_4$ & $a^2$ & ${\tt -0.2874798361}$   \\ \hline
$ A_5$ & $a^{5/2}$ & ${\tt  -22.04012159}$  \\ \hline 
$B_3$ & $a^3\log a$ &${\tt -6.966574583 }$ \\ \hline
$ A_6$ & $a^3$ & ${\tt -6.250578833 }$ \\ \hline
$ A_7$ & $a^{7/2}$ &${\tt  39.27382902 }$  \\ \hline
$B_4$ & $a^4\log a$ &${\tt 10.92891913  }$ \\ \hline
$ A_8$ & $a^4$ &${\tt  19.31687979 }$  \\ \hline
$ A_9$ & $a^{9/2}$ & ${\tt -82.17477248 }$  \\ \hline
$ B_5$ & $a^5\log a$ &${\tt -20.01663759 }$  \\ \hline
$ A_{10}$ & $a^5$ & ${\tt -50.35178499 }$ \\ \hline
$A_{11}$ & $a^{11/2}$ & ${\tt 186.9039750  }$\\ \hline
$B_6$ & $a^6\log a$ &${\tt 40.63376347  }$\\ \hline
$ A_{10}$ & $a^6$ &${\tt  128.7996196 }$\\ \hline
$ A_{11}$ & $a^{13/2}$ &${\tt -444.7385574 }$\\ \hline
\end{tabular}
\begin{tabular}{|c|c|c|}
\hline
 coef. & scale & value \\ \hline
$ M_0$ &  $ 1 $ & ${\tt -13.35083105  }$   \\ \hline
$M_1$ &  $\mu^{1/4} $  & ${\tt -35.34533603 }$    \\ \hline
$N_1$ & $ \mu^{1/2}\log \mu  $  &${\tt -9.183673469 }$   \\ \hline
$M_2$ & $ \mu^{1/2}  $ & ${\tt -25.71572403  }$ \\ \hline
$M_3$ & $ \mu^{3/4} $ & ${\tt 60.69366755  }$\\ \hline
$ N_2$ & $\mu \log \mu  $ & ${\tt -41.92449575 }$   \\ \hline
$ M_4$ & $ \mu$ & ${\tt -215.4221683}$   \\ \hline
$ M_5$ & $ \mu^{5/4}$ & ${\tt -45.92851439 }$ \\ \hline
$ N_3$ & $ \mu^{3/2}\log \mu$ &${\tt - 242.5333437}$ \\ \hline
$M_6$ & $  \mu^{3/2} $ & ${\tt - 960.8699623}$ \\ \hline
$M_7$ & $ \mu^{7/4} $ & ${\tt 755.3601690 }$ \\ \hline
$N_4$ & $ \mu^{2}\log \mu$ &${\tt -1587.303140}$ \\ \hline
$M_8$ & $\mu^2$ & ${\tt -3308.441120}$ \\ \hline
$M_9$ & $\mu^{9/4}$ & ${\tt 1090.837521  }$  \\ \hline
$N_5$ & $\mu^{5/2}\log \mu$ &${\tt -11017.80445}$  \\ \hline
$ M_{10}$ & $ \mu^{5/2}$ & ${\tt -134120.3771}$ \\ \hline
$ M_{11}$ & $\mu^{11/4}$ & ${\tt  22519.75418 }$\\ \hline
$N_6$ & $ \mu^3\log \mu $ & ${\tt  -79363.78673 }$\\ \hline
$ M_{10}$ & $ \mu^3 $ & ${\tt   904656.6104}$ \\ \hline
$M_{11}$ & $\mu^{13/4}$ &${\tt 87833.05069}$\\ \hline
\end{tabular}
\caption{The 
20 first coefficients of the asymptotic expansion
for the Bogdanov map (left, ${\tilde \gamma}=3$) and
the Quadratic map (right, ${\gamma}=-3$).
All the given digits are correct.
}
\end{table}

\subsection{H\'enon map}
The last example to be considered in this paper
is 
the H\'enon map \cite{henon} 
defined by
\begin{eqnarray*}
{\bf H}={\bf H}_{{\tilde a},{\tilde b}} & : 
& {\bR}^2 \rightarrow {\bR}^2, \ (u,v) \mapsto (u_1,v_1)  
\end{eqnarray*}  
where 
$$
u_1 = v, \ \ v_1 = {\tilde a}v^2-{\tilde b} u+1.
$$
See \cite{ks} for recent results concerning this family.
The H\'enon map has 
a fixed point of Bogdanov Takens type at 
${\tilde a}={\tilde b} =1$.
We chose ${\tilde a}$ 
as the
main parameter and ${\tilde b}$ 
as the slave parameter. We note that
the H\'enon map is conjugate to the Bogdanov family
in the special case of ${\tilde\gamma}=0$. 
The conjugacy is given by 
the following change of coordinates and parameters  
$$
u=x, 
\ \ v=x+y+x^2 +a x+by, 
\ \ \ {\tilde b}=b+1,\ {\tilde a}= (1+b/2)^2-a^2/4.  
$$
We also observe that 
the H\'enon map can be transformed 
to the form (\ref{unf}) with the non linear term 
of the form (\ref{unfold})
 by putting 
$$
u=
\frac{1}{\tilde a}(x +\frac{{\tilde b}+1}{2}),\  
v= 
\frac{1}{\tilde a}(x +\frac{{\tilde b}+1}{2}) +  \frac{1}{\tilde a} y.
$$
In the new system of coordinates,
the H\'enon map takes the  form (\ref{unfold}) with
$$
f_{\mu,\nu} =(x+y)^2-\mu +\nu y,
$$
where 
$$
\mu =(1+\frac{\nu}{2})^2-{\tilde a}, 
\ \ \nu = {\tilde b}-1.
$$
The H\'enon map admits a homoclinic zone near the line
$$
{\tilde b}^\pm({\tilde a}) \equiv 1, \ {\tilde a} \geq 1.
$$
In the case of the H\'enon map we define the normalized width 
of the zone by
\begin{eqnarray*}
 {\tilde S} ({\tilde a})
& = &  
\frac{{\tilde b}^+({\tilde a})-{\tilde b}^-({\tilde a})}{K(1-{\tilde a}, 0)}.
\end{eqnarray*}
Our numerical experiments show that 
$ {\tilde S}$ 
has the following
asymptotic expansion:
\begin{eqnarray}
\label{bog-diverge}
{\tilde S} ({\tilde a}) & = & 
\sum_{k\geq 0}{\tilde A}_k (1-{\tilde a})^{k/4}.
\end{eqnarray}
 
\noindent
Unlike the case of the Bogdanov map with $\gamma {\not =2}$ (i.e., ${\tilde \gamma} {\not =0}$),
the asymptotic expansion does not contain logarithmic terms.
We expect
this property to be closely related  
to the fact that the H\'enon map contains a one parametric subfamily
of area preserving maps.
In general, even when $\gamma=2$, there is no reason to expect the logaritmic terms to vanish
for a map $F_{\mu,\nu}$.
\begin{table}
\begin{center}
\begin{tabular}{|c|c|c|}
\hline  
 coef. & scale & H\'enon map 
\\ \hline
${\tilde A}_0$ &  $ 1 $ &        $2.4744255935532510538408*10^6$\\ \hline
${\tilde A}_1$ &  $|{\tilde a}-1|^{1/4} $  &  ${\tt -2.878113364919828141704 *10^6}$\\ \hline
$ {\tilde A}_2$ & $  |{\tilde a}-1|^{1/2} $  &    ${\tt 1.8211174314566012763528*10^6}$ \\ \hline
${\tilde A}_3$ & $ |{\tilde a}-1|^{3/4} $ & ${\tt -412552.07921345800366019\ \ \  \ \ \ \  \ \  }$ \\ \hline
${\tilde A}_4$ & $ |{\tilde a}-1|$ &       $ {\tt -309961.28583121907079391\ \ \  \ \ \ \  \ \ }$\\ \hline
${\tilde A}_5$ & $|{\tilde a}-1|^{5/4}$ &   $ {\tt 257055.93487794037812901\ \ \  \ \ \ \  \ \ }$ \\ \hline
${\tilde A}_6$ & $|{\tilde a}-1|^{3/2}$ &  ${\tt -56830.201956139947433580 \ \ \  \ \ \ \  \ \ }$ \\ \hline
${\tilde A}_7$ & $|{\tilde a}-1|^{7/4}$ &  ${\tt-12386.990577003086404843 \ \ \  \ \ \ \  \ \ }$\\ \hline 
${\tilde A}_8$ & $|{\tilde a}-1|^{2}$ &    ${\tt-11792.964908478734939516  \ \ \  \ \ \ \  \ \ }$\\ \hline
${\tilde A}_9$ & $|{\tilde a}-1|^{9/4}$ &  $ {\tt 18742.189161591275288347  \ \ \  \ \ \ \  \ \ }$\\ \hline
${\tilde A}_{10}$ & $|{\tilde a}-1|^{5/2}$  & ${\tt-4774.6727458595190485600 \ \ \  \ \ \ \  \ \ }$  \\ \hline
${\tilde A}_{11}$ & $|{\tilde a}-1|^{11/4}$ & ${\tt-2822.9663193640187675835 \ \ \  \ \ \ \  \ \ }$\\ \hline
${\tilde A}_{12}$ & $|{\tilde a}-1|^{3}$ &    $ {\tt3276.6438736125169964394 \ \ \  \ \ \ \  \ \ } $\\ \hline
${\tilde A}_{13}$ & $|{\tilde a}-1|^{13/4}$ & ${\tt-1910.5466958542171966392 \ \ \  \ \ \ \  \ \ }$\\ \hline
${\tilde A}_{14}$ & $|{\tilde a}-1|^{7/2}$ &  ${\tt 7704.6605615546853854041\ \ \  \ \ \ \  \ \ }$\\ \hline
${\tilde A}_{15}$ & $|{\tilde a}-1|^{15/4}$ & ${\tt-7827.0351891507566506398 \ \ \  \ \ \ \  \ \ }$\\ \hline
${\tilde A}_{16}$ & $|{\tilde a}-1|^{4}$ &    ${\tt 13919.102717097324631620\ \ \  \ \ \ \  \ \ }$ \\ \hline
${\tilde A}_{17}$ & $|{\tilde a}-1|^{17/4}$ & ${\tt-11932.139780641352182621\ \ \  \ \ \ \  \ \ }$\\ \hline
${\tilde A}_{18}$ & $|{\tilde a}-1|^{9/2}$ &  ${\tt 22120.721696311178434645 \ \ \  \ \ \ \  \ \  }$\\ \hline
\end{tabular}
\end{center}
\caption{The $19$ first coefficients 
in (\ref{bog-diverge-henon}).
All the given digits are correct.
We also conjecture that the series (\ref{bog-diverge-henon}) belongs to the Gevrey-1 class.
}
\end{table}

\section{Computing the width of the homoclinic zone}
In this section, our approach concerns
the   Quadratic family 
${\bf Q}_{\mu,\nu,\gamma}$. 
The other
families (Bogdanov and H\'enon) are treated 
in a similar way.
From now on, we do not mention the $(\mu,\nu,\gamma)$ 
dependences when it is not necessary, but we may
emphasise that dependence when it is needed.

\subsection{Strategy}

\begin{itemize}
\item[i)] We assume an ansatz and in particular the one given in 
formula (\ref{teo-log});
\item[ii)] Compute ${\tilde n}$ (several hundreds) values of the 
width for values of $\mu^{1/4} \in [c, d]$,
where $0<c<d$ are close to $0$  
(typically $c \approx 5/1000$, $d \approx 1/100$). 
It is convenient to work 
with the so called `normalised 
width of homoclinic zone' defined by
\begin{eqnarray*}
S_\gamma({\mu}) & =& 
\frac{\nu^+(\mu) -\nu^-(\mu)}{K(\mu, \gamma-2)}
\end{eqnarray*}
where
$K$ is defined by (\ref{kmugamma}). 
The result is collected in a set of data of the form
\begin{eqnarray}\label{zz-h}
{ {\cal H}} =
\{
(\mu_i^{\frac{1}{4}},\log(S_\gamma({\mu_i})), \ c\leq \mu_i \leq 
 d, \ i=1,\ldots,{\tilde n}
\}.
\end{eqnarray}

\item[iii)] 
Take $\ell \in {\bN}$ such that $3\ell/2+1 \leq {\tilde n}$ and $\ell >>1$ even.
Then we compute the coefficients $M_k$, $k=0,\ldots,\ell$ and
$N_k$, $k=1,\dots, \ell/2$,
of the truncated expansion
$$
G^{\{{ 3\ell/2}\}} (\mu) 
=\sum_{k=0}^\ell M_k(\gamma)\mu^{k/4}+
\log \mu \sum_{k=1}^{\ell/2} N_k(\gamma)\mu^{k/2}
$$
to 
interpolate the set ${\cal H}$, i.e., 
for all integer $i=1,\ldots, 3\ell/2 +1$, we have
$$
\log S_\gamma({\mu_i}) =
\sum_{k=0}^\ell M_k(\gamma)\mu_i^{k/4}+\log \mu_i 
\sum_{k=1}^{\ell/2} N_k(\gamma)\mu_i^{k/2}.
$$
See subsection 3.10 for more details.

\end{itemize}
\vskip10pt
\noindent
{\bf Remarks}:
\begin{itemize}
\item
For the Bogdanov family, the set of data for the normalised width is
denoted by
\begin{eqnarray}\label{bog-nom}
{\tilde {\cal H}} & = & \{ (a_i^{\frac{1}{2}}, \log( {\tilde S}_\gamma(a_i) )), 
\ {\tilde c}<a_i< {\tilde d}, \ 
 \ i=1,\ldots,{\tilde n}\},
\end{eqnarray}
where $$ {\tilde S}_\gamma(a_i)
 = \frac{b^+(a_i)-b^-(a_i)}{K(a_i^2/4,{\tilde \gamma})}  \ {\mathrm{and}}\  0<{\tilde c}<{\tilde d}.
$$
\item
For the H\'enon family, the set of data for the normalised width is
denoted by
\begin{eqnarray}\label{t-c-z}
{\tilde {\cal Z}} & = & \{ (|1-{\tilde a}_i|, {\tilde S}({\tilde a}_i)), 
\ {\tilde c}<|1-{\tilde a}_i|< {\tilde d}, \ 
 \ i=1,\ldots,{\tilde n}\},
\end{eqnarray}
where 
$$
 {\tilde S}({\tilde a}_i) = 
\frac{{\tilde b}^+({\tilde a}_i)-{\tilde b}^-({\tilde a}_i)}{K((1-{\tilde a}_i),0)}  
\ {\mathrm{and}}\  0<{\tilde c}<{\tilde d}.
$$
\end{itemize}

    \subsection{Invariant manifolds}
We now 
compute the stable and unstable manifold at the saddle point. 
In what follows, our description concerns the
Quadratic map ${\bf Q}$ 
but similar computations are done for the Bogdanov map and the H\'enon map.

From (\ref{rfam}), the map ${\bf Q}$ has two
fixed points 
$${\bf S}_\mu=(\sqrt{\mu},0), \ {\mathrm{and}} \ {\bf C}_\mu=(-\sqrt{\mu},0).$$ 
${\bf C}_\mu$ is a focus and 
${\bf S}_\mu$ is a saddle and will be
the point of interest.
The eigenvalues of $d{\bf Q}({\bf S}_\mu)$
are given by
\begin{eqnarray*}
\lambda_1 & = & \frac{1}{2}
\biggl(2+\nu+\gamma \sqrt{\mu}-{\sqrt{(\gamma {\sqrt{\mu}}+\nu)^2 +8{\sqrt{\mu}}}}\biggr), \\
\lambda_2 
& = & 
 \frac{1}{2}\biggl(2+\nu+\gamma\sqrt{\mu}
+
{\sqrt{(\gamma {\sqrt{\mu}}+\nu)^2 +8{\sqrt{\mu}}
}}\biggr).  
\end{eqnarray*}
For $\mu>0$ sufficiently small it is clear that $\lambda_1<1<\lambda_2$. 
At the saddle ${\bf S}_\mu$,  the Taylor expansion 
of the local  stable manifold $W^s_{\mathrm{loc}}$
 and that of the local unstable manifold
$W^u_{\mathrm{loc}}$
are computed as follows. Denote by 
$$\Phi_s  : \ ({\R},0) 
\rightarrow ({\R}^2, {\bf S}_\mu),
\ z \mapsto \Phi_s(z)= \biggl( \sqrt{\mu}+\sum_{k=1}^{\infty} \varphi_k z^k,
\sum_{k=1}^{\infty} \psi_k z^k
\biggr)$$
$$\Phi_u  : \ ({\R},0) \rightarrow ({\C}^2,{\bf S}_\mu), \ z \mapsto \Phi_u(z)= 
\biggl( \sqrt{\mu}+\sum_{k=1}^{\infty} f_k z^k,
\sum_{k=1}^{\infty} p_k z^k
\biggr)$$
the parameterisations which respectively satisfy 
\begin{eqnarray}\label{satie}
\Phi_s( \lambda_1 z)  = {\bf Q} \circ \Phi_s(z) & {\mathrm {and}} \ &
    \Phi_u( \lambda_2 z)  =  {\bf Q} \circ \Phi_u(z)
\end{eqnarray}
for all $z$ near $0$. Substituting the series into
 (\ref{satie}) and collecting  terms of the same order in $z$ 
 we get
\begin{eqnarray}\label{vstable}
\left\{
\begin{array}{rcl}
\varphi_k +\psi_k & = & \lambda_1^k \varphi_k, \ \ \ k\geq 1 \\
\sum_{j=0}^k  \varphi_j \varphi_{k-j} +\gamma \sum_{j=0}^k  \varphi_j \psi_{k-j}
+\nu \psi_k
& = &  \lambda_1^k \psi_k
\end{array}
\right. 
\\ \label{vunstable}
\left\{
\begin{array}{rcl}
p_k +f_k & = & \lambda_2^k p_k, \ \ \ \ \ \ k\geq 1 \\
\sum_{j=0}^k  f_j f_{k-j} +\gamma \sum_{j=0}^k  f_j p_{k-j}
+\nu p_k
& = &  \lambda_2^k p_k
\end{array}
\right.
\end{eqnarray}
Since $\lambda_2>1$ and ${\bf Q}$ is entire, 
from (\ref{satie}) we easily deduce that the
radius of convergence of the
series defined in (\ref{vunstable}) 
is infinite.
Denote by $\rho$ the radius of convergence
of the series defined in (\ref{vstable}).
We fix $N_{\mathrm {max}} \in {\N}$. 
Since we are after a single branch of the stable  
manifold 
we write
\begin{eqnarray*}
W^s_{\mathrm{loc}}  \approx W^s_{N_{\mathrm{max}}}
=\{\Phi_{s,N_{\mathrm{max}}}(z), \ 0 \leq z \leq \delta_s \},
\end{eqnarray*}  
where
$$
\Phi_{s,N_{\mathrm{max}}}(z)=\biggl( \sqrt{\mu}+\sum_{k=1}^{N_{\mathrm{max}}} \varphi_k z^k,
\sum_{k=1}^{N_{\mathrm {max}}} \psi_k z^k
\biggr)
$$  
and where
$0<\delta_s<\rho$. 
We proceed  in the same way for the local unstable manifold, i.e., 
\begin{eqnarray*}
W^u_{\mathrm {loc}}  \approx W^u_{N_{{\mathrm {max}}}}=\{ 
\Phi_{u,N_{\mathrm{max}}}(z),
\ 0 \leq z \leq \delta_u \}.
\end{eqnarray*}  
where 
\begin{eqnarray}
\label{ust}
\Phi_{u,N_{\mathrm{max}}}(z)= \biggl(\sqrt{\mu}+\sum_{k=1}^{N_{\mathrm {max}}} f_k z^k,
\sum_{k=1}^{N_{\mathrm {max}}} p_k z^k
\biggr)
\end{eqnarray}
and where
$0<\delta_u<<1$. 
The local invariant 
manifolds are computed with the following precision:
\begin{eqnarray*}
\|\Phi_{s,N_{\mathrm{max}}}(z)-\Phi_{s}(z)\| 
 =   {\cal O}
(z^{N_{\mathrm{max}}}), \qquad 
\| \Phi_{u,N_{\mathrm{max}}}(z)-\Phi_{u}(z) \| 
 = {\cal O}(
 z^{N_{\mathrm{max}}}
).
\end{eqnarray*}
In particular we have 
\begin{eqnarray}\label{py0}
\| \Phi_u(\lambda_2 z)-{\bf Q} \circ \Phi_{u, N_{\mathrm{max}}}(z) \| ={\cal O}(z^{N_{\mathrm{max}}}).
\end{eqnarray} 
Since we need to study the map 
when
homoclinic orbits are present, 
we need a good estimate of 
the 
global unstable manifold.
Recall that $\Phi_u$ is entire and therefore both components defined in 
(\ref{ust}) converge for all $z$ as $N_{\mathrm{max}} \to \infty$. 
However, for large $z$, 
the computation of the unstable manifold requires too many
coefficients and therefore
(\ref{ust}) is not very 
convenient. 
We then proceed as follows.
%
Let ${\mathrm {P}_0}= \Phi_u(z_0) \in W^u$ and choose $m_0$
such that 
$$
z_1=\lambda_2^{-m_0}z_0 \leq \delta_u.
$$  
Then, for any fixed $m_0$, we have
\begin{eqnarray*}
P_0 & = & \lim_{N_{\mathrm{max}} \to \infty} 
{\bf Q}^{m_0} \circ \Phi_{u,N_{\mathrm{max}}}( z_1) 
\end{eqnarray*}
and if $z_1 <<1$  the convergence 
is fast.
%
Therefore, by putting 
\begin{eqnarray*}
W^u \approx W^u_{N_{\mathrm{max}},m}= \{ 
{\bf Q}^m \circ \Phi_{u,N_{\mathrm{max}}}(\lambda_2^{-m} z), 
\ 0 \leq   z \leq { z_0}\}, \  m \geq { m_0}
\end{eqnarray*} 
we get an accurate  estimation 
of the global unstable manifold.

   \subsection{Jacobian and Wronskian functions}
Before introducing the splitting function 
which will play a key role in the paper, we need to introduce two
additional functions.
We first define
\begin{eqnarray*}\label{jacobsim}
J & : & {\cal D} \rightarrow \C, \ z \mapsto \det d{\bf Q} (\Phi_s(z))
\end{eqnarray*}
as the Jacobian of the map ${\bf Q}$
along the stable manifold 
$$\Phi_s(z)=(\Phi_{s,x}(z),\Phi_{s,y}(z)).$$
A straightforward computation gives
\begin{eqnarray}\label{jacob}
J(z)=1+\nu +(\gamma-2) \Phi_{s,x} (z)-\gamma  \Phi_{s,y} (z).
\end{eqnarray}
In terms of series, from (\ref{jacob}) we get
\begin{eqnarray}\label{j-serie}
J(z)  =   \sum_{k=0}^{\infty} J_k z^k, & \ {\mathrm{where}} & \ \\ \nonumber
J_0  =   1+\nu +(\gamma-2) \sqrt{\mu}, & \ {\mathrm {and}} \ &  
\forall k>0,\  
J_k=(\gamma-2)\phi_k -\gamma \psi_k . 
\end{eqnarray}
The Wronskian function (along the local stable manifold)
\begin{eqnarray*}\label{jacobsim2}
\Omega & : & {\cal D} \rightarrow \R, \ z \mapsto \Omega(z)
\end{eqnarray*}
satisfies
\begin{eqnarray}\label{w-pro}
\Omega(\lambda_1 z) & = & J(z)\Omega(z).
\end{eqnarray}
We put $\Omega_0=1$ and look for a solution of (\ref{w-pro})
of the form
\begin{eqnarray}\label{w-form}
\Omega(z) = z^{{\log J_0}/{\log\lambda_1}}
\biggl(1 + \sum_{k=1}^{\infty}\Omega_k z^k 
\biggr).
\end{eqnarray}
With (\ref{w-pro}), (\ref{j-serie}),
and (\ref{w-form}), it follows that
$$
\Omega_n
=
\frac{1}{\lambda_1-J_0}
\biggl(
 J_n+\sum_{j=0}^{n-1}\Omega_j J_{n-1-j} 
\biggr).
$$
Both series (\ref{j-serie}) and (\ref{w-form}) are convergent.
The functions $J$ and $\Omega$ will be approximated by
\begin{eqnarray*}
J_{N_{\mathrm{max}}}(z) 
= \sum_{k=0}^{N_{\mathrm{max}}} J_k z^k\ 
\ {\mathrm {and}} \ \ \ \Omega_{N_{\mathrm{max}}}(z) = z^{{\log J_0}/{\log\lambda_1}}
\biggl(1 + \sum_{k=1}^{N_{\mathrm{max}}}\Omega_k z^k
\biggr)
\end{eqnarray*}
respectively. In this way, 
we have
 \begin{eqnarray}
\label{py2}
|\Omega_{N_{\mathrm{max}}}(\lambda_1 z)-
J_{N_{\mathrm{max}}}(z)\Omega_{N_{\mathrm{max}}}(z)
| & = & {\cal O}(|z|^{N_{\mathrm{max}}}). 
\end{eqnarray}
 
    \subsection{Splitting function and flow box theorem}
In this section, we introduce the key part of our 
techniques. Recall that in our investigation for the width of the homoclinic zone, we fix the 
value of the main parameter and look for values $\nu^+$ and $\nu^-$ of the slave parameter 
that correspond, respectively, to the first and the last homoclinic tangency.
In order to find a homoclinic point we need
to adjust the slave parameter in such a way that two curves on the plane
have an intersection. Finding a homoclinic tangency requires additional
adjustments to make this intersection degenerate. This problem is much
easier in the discrete flow box coordinates, in which the stable
curve coincides with the horizontal axis and the unstable one is a graph
of a periodic function. A further simplification will be achieved
by observing that this periodic function is very close to a trigonometric
polynomial of the first order.
The splitting function $\Theta=\Theta_{\mu, \nu}$ 
we shall introduce now is such that 
 the first and the last tangency correspond to double zeroes
of $\Theta_{\mu, \nu^+}$ and $\Theta_{\mu, \nu^-}$ respectively.
Our investigation amounts then to
 finding values $\nu^+$ and $\nu^-$ such that
$\Theta_{\mu, \nu^+}$ and $\Theta_{\mu, \nu^-}$ possess double zeroes. 

In this section, we present the splitting function $\Theta_{\mu,\nu}$
for the Quadratic map, in the case of the Bogdanov map, the splitting 
function is denoted by
$\Theta_{a,b}$.
In what follows, we 
assume that the parameter $(\mu, \nu)$ 
is such that the map  ${\bf Q}$
possesses a homoclinic orbit, i.e., 
the unstable manifold intersects the local stable manifold at
a point $\Phi_u(z_u)=q_0=\Phi_s(z_s)$.
Then we fix  a neighbourhood ${\cal U}$ of the point $q_0$.
 We parametrise
$W^s_{\mathrm {loc}}$ near $q_0$ 
by
\begin{eqnarray*}
{\bf \Gamma}_s : {\bI}_0 & \mapsto  & \Phi_s(z_s\cdot \lambda_1^t) 
\end{eqnarray*}
where ${\bI}_0=(-1, 1)$
and 
$W^u$ near $q_0$ by
\begin{eqnarray*}
{\bf \Gamma}_u : {\bI}_0 & \mapsto  & \Phi_u(z_u\cdot \lambda_2^t).
\end{eqnarray*}
 Now we state the following
 (flow box)
lemma \cite{Gelfreich1996}.
\begin{lemma}
{\rm
There exists ${E}_0>0$  
and an analytic diffeomorphism
\begin{eqnarray*}
\Psi  :  (-{E}_0, {E}_0) \times  {\bI}_0 & \rightarrow& {\bR}^2, \ 
 \\ 
(E,t) &\mapsto & \Psi(E,t) = (X(E,t), Y(E,t))
\end{eqnarray*}
such that the following hold
\begin{itemize}
          \item[\rm{i)}] $\Psi(E,t+1) = {\bf Q} \circ \Psi(E,t)$,
             \item[\rm{ii)}] $\Psi(0,0)=q_0$, 
$\Psi(0,t) \in W^s_{\mathrm{loc}}$ for  $t\in  {\bI}_0$,
          \item[\rm{iii)}] the Jacobian matrix
\begin{eqnarray}\label{JJ}
d\Psi(E,t)=
\displaystyle\left(
\begin{array}{cc}
{\partial X }/{\partial E} &{\partial X }/{\partial t} \\
{\partial Y }/{\partial E} &{\partial Y }/{\partial t}
\end{array}\right)
,
\end{eqnarray}
is such that the second column of 
$ d\Psi(0,t)$
 is 
$\dot {\bf \Gamma}_{s}=d{\bf \Gamma}_s(t)/dt$, 
          
          \item[\rm{iv)}]  the map ${\hat \Omega}(E,t)=\det d\Psi(E,t)$ 
satisfies
${\hat \Omega}(0,t)=\Omega(z_s\cdot \lambda_1^t)$;
\end{itemize}
}
\end{lemma}

\vskip10pt
\noindent
The splitting function, denoted by $\Theta_{\mu,\nu}(t)$,
 is the first component of 
$$
 \Psi^{-1} \circ {\bf \Gamma}_u(t) -\Psi^{-1} \circ {\bf \Gamma}_s(t).
$$
Applying Taylor  theorem at the stable manifold, we
get 
\begin{eqnarray}
\nonumber
\Psi^{-1} \circ {\bf \Gamma}_u(t) -\Psi^{-1} \circ {\bf \Gamma}_s(t) 
 & =&   
d\Psi^{-1}(\Psi(0,t))\cdot
\biggl(
{\bf \Gamma}_u(t)-{\bf \Gamma}_s(t)
\biggr)\\ 
\label{JJ2}
& + & 
{\cal O}
\biggl(
\|{\bf \Gamma}_u(t)-{\bf \Gamma}_s(t)\|^2
\biggr).
\end{eqnarray}
 The following properties hold: 
\par\vskip10pt
[-] Let  $0<{\tilde \delta}<\pi$. The map 
 $\Theta_{\mu,\nu}$ has an analytic continuation 
onto the rectangle:
\begin{eqnarray}\label{domain-b}
{\mathrm {B}} = 
\{ t \in {\C} 
\ | \ t=t^\prime+i t^{\prime\prime}, \ t^\prime \in  {\bI}_0 , 
\ |t^{\prime\prime}| \leq \varrho\}, 
\ |\varrho| < (\pi-{\tilde \delta})/ |\log \lambda_1|.
\end{eqnarray}
The function $\Theta_{\mu,\nu}$ is periodic so we can expand it into Fourier
series:
 \begin{eqnarray*}    
{ \Theta}_{\mu,\nu}(t)  =  
\sum_{j=-\infty}^{\infty} {{\gmp}}_{j}(\mu,\nu) e^{2i \pi t}.  
\end{eqnarray*}
As usual, the Fourier coefficients are defined by an integral:
\begin{eqnarray*}
{{\gmp}}_k(\mu,\nu) = \int_0^1 
{ \Theta}_{\mu,\nu}(t)e^{-2i k\pi t}dt, 
\ \  {\mathrm { \ for \ each }} \ k \in {\N}.
\end{eqnarray*}
Let $0<\rho<(\pi-{\tilde \delta})/|\log(\lambda_1)|$. 
Since the integral of 
$
{ \Theta}_{\mu,\nu}(t)e^{-2ik\pi t}
$ 
over the boundary of the rectangle
$
\{(t^\prime +i t^{\prime\prime}) 
\ | \ 0\leq t^\prime\leq 1,\  0\leq  t^{\prime\prime} \leq \rho \}
$
vanishes, we conclude 
 \begin{eqnarray}
\label{f_3}
 \int_0^1 { \Theta}_{\mu,\nu}(t)e^{-2ik\pi t}dt & = &  
 e^{-2k \pi \rho} \int_0^1 { \Theta}_{\mu,\nu}(t+i\rho)e^{-2i k\pi t}dt.
\end{eqnarray}
Consequently  
\begin{eqnarray}\label{h-k}
|{{\gmp}}_k(\mu,\nu) | \leq \sup_{t \in  {\bI}_0} 
|\Theta_{\mu,\nu} (t+ i\rho)|\cdot e^{-2 |k|\pi\rho}, 
\end{eqnarray}
i.e., the harmonics of $\Theta_{\mu,\nu}$
decrease exponentially. 
The function $\Theta_{\mu,\nu}$
can be well approximated by the sum of zero and first order harmonics:
 \begin{eqnarray}
\label{w-w}
\Theta_{\mu,\nu}(t) & =&  
{{\gmp}}_{-1}(\mu,\nu) e^{-2i\pi t} 
+  
{{\gmp}}_{0}(\mu,\nu) 
+ 
{{\gmp}}_{1}(\mu,\nu) e^{2i\pi t}
 +  {\cal O}_2(t)
\end{eqnarray}
or equivalently,  $\Theta_{\mu,\nu}$
is well approximated 
by a trigonometric polynomial function 
\begin{eqnarray}\label{real-exp}
\Theta_{\mu,\nu}(t)  =  
{{\gmp}}_{0}(\mu,\nu) +2 |{{\gmp}}_{-1}(\mu,\nu)|
\cos(2\pi t +{\mathrm{arg}}({{\gmp}}_{-1}(\mu,\nu))) 
 + {\cal O}_2(t)
\end{eqnarray}
where
\begin{eqnarray}
\label{o-2}
 \sup_{i\in {\bI}_0}
|{\cal O}_2|(t)
= {\cal O}(\sup_{t\in  {\bI}_0}|\Theta_{\mu,\nu}(t)|^2).
\end{eqnarray}
\par
\vskip10pt
[-] Since 
$d\Psi^{-1}(\Psi(0,t))
=(d\Psi(0,t))^{-1},$ 
we have
 \begin{eqnarray*}
d\Psi^{-1}(\Psi(0,t))
& = & \frac{1}{{\hat {\Omega}}(0,t)}
\displaystyle\left(
\begin{array}{cc}
{\partial Y }/{\partial t} & -{\partial X }/{\partial t} \\
{\partial X }/{\partial E} & {\partial Y }/{\partial E}
\end{array}\right).
\end{eqnarray*}
Furthermore,  
$$
\Psi^{-1}(\Gamma^u(t))=
\biggl(E_u(t), T_u(t) \biggr), \  
\Psi^{-1}(\Gamma^s(t))
=\biggl(E_s(t), T_s(t) \biggr)=(0,t),
 $$ 
with (\ref{JJ}) and (\ref{JJ2}) it follows that
\begin{eqnarray}\nonumber
\Theta_{\mu,\nu}(t)=E_u(t)-E_s(t) & =& \frac{1}{{\hat {\Omega}}(0,t)}
\det \left(\begin{array}{ccc}
\displaystyle\frac{d}{dt} 	{\bf \Gamma}_s(t) & , & {\bf \Gamma}_u(t) -{\bf  \Gamma}_s(t) 
 \end{array}\right) \\  \label{def-tet}
& +& 
{\cal O}(\|{\bf \Gamma}_u(t) -{\bf  \Gamma}_s(t) \|^2).
\end{eqnarray}
Thus, 
we obtain a formula suitable for computation of the splitting function
in terms of the parametrization of the stable and unstable manifold:
%
\begin{eqnarray}
\label{formula}
\Theta_{\mu,\nu}(t) = {\tilde \Theta}_{\mu,\nu}(t) + {\tilde h}_{\mu,\nu}(t)
\end{eqnarray}
where
\begin{eqnarray}
\label{for-tilde}
 {\tilde \Theta}_{\mu,\nu}(t) 
& = & \displaystyle\frac{1}{{{\Omega}}(z_s\cdot\lambda_1^t)} 
\det \left(\begin{array}{ccc}
\displaystyle\frac{d}{dt}    {\bf \Gamma}_s(t) &  
{\bf \Gamma}_u(t) -{\bf  \Gamma}_s(t) 
 \end{array}\right)
\end{eqnarray}
is the splitting determinant 
and 
\begin{eqnarray}
\label{til-h-split}
 |{\tilde h}_{\mu,\nu}(t)| & = & {\cal O}(\sup_{t \in {\bI}_0}|
{\tilde \Theta}_{\mu,\nu}(t) |^2).
\end{eqnarray}
Note that 
 even if the invariant manifolds and the Wronskian are computed 
with a very high precision,
the function $\Theta_{\mu,\nu}(t)$ 
is only evaluated 
with  a relative error of 
order  ${\cal O}(\sup_{t\in {\bI}_0}|\Theta_{\mu,\nu}|)$. 
%

\subsection{Approaching a primary homoclinic orbit}
In order to compute the width of the homoclinic zone, we first find 
a value $\nu={\bar \nu}$ where the map possesses a primary homoclinic orbit. 
Near $\nu ={\bar \nu}$, Lemma 1 will then be  applied and the
splitting determinant
${\tilde \Theta}_{\mu,\nu}$ will be computed. We proceed as follows:
we fix $0<z_s<\delta_s$  and a section 
$\Sigma$ transverse to the local stable manifold at  $p_\nu = \Phi_s(z_s)$.
We parametrise $\Sigma$ as follows
$$
\Sigma= \{p_\nu + (0,y), \ -y_0 <y< y_0  \}
$$
where $0<y_0<<1$.
For each value of the main parameter,
we consider the slave parameter being close to 
$\nu_0 =(5(\gamma-2)/7){\sqrt{\mu}}$ and compute 
a point $q_\nu \in W^u\cap \Sigma$ which is 
the `first intersection'
of 
 $W^u$ with the section. 
In order to increase the speed 
of computations we use Newton's method to solve the equation 
${\bf \Gamma}_u(t) \in \Sigma$. After that we adjust $\nu$
in such a way that $q_\nu=p_\nu$.
We do not know an easy way to evaluate the derivative of $q_\nu$
with respect to $\nu$, therefore we cannot apply Newton's method.
However, we replace the derivative by a finite difference approximation and use 
the so called `secant' method. In other words we consider
the limit of the following sequence:
$$
\nu_{n+1}= \nu_{n} +
\frac{{\bar \delta}y_{\nu_n}}{y_{\nu_n+{\bar \delta}}-y_{\nu_n}}
$$
where $q_\nu=p_\nu+(0,y_\nu)$ and 
where $0<{\bar \delta}<<1$.
Denote by $${\bar  \nu}= \lim_{n\to \infty} \nu_n.$$
Since $p_{\bar \nu} =q_{\bar \nu}$, the point $(\mu, {\bar \nu})$
belongs to the homoclinic zone. 

\vskip10pt
Our next step is with the computation of the width $\nu^+(\mu)-\nu^-(\mu)$
for the given value of $\mu$. 
The zeroes (and double zeroes) 
of $\Theta_{\mu,\nu}$ 
are in one to  one correspondence with  primary homoclinic orbits 
(and homoclinic tangencies)
for 
the corresponding map, 
see \cite{Gelfreich1996,Gelfreich2003} for more details.
We then replace the problem of 
finding homoclinic points and homoclinic tangencies
by finding 
double zeroes of the splitting function $\Theta_{\mu,\nu}$.
%
%
%
%
%
\subsection{First and last tangency}
The most natural way to compute
the width of homoclinic zone is to estimate both $\nu^+=\nu^+(\mu)$ and $\nu^-=\nu^-(\mu)$. 
Write 
\begin{eqnarray}
\label{hat-theta}
{\Theta}_{\mu,\nu}(t) & = &  {\bf P}_0(\mu,\nu) + {\hat {\Theta}}_{\mu,\nu}(t).
\end{eqnarray}
At the first tangency, ($\nu=\nu^-$) the graph of the splitting function 
is located below the $t$ axis and  $\Theta_{\mu,\nu^-}$ 
admits  a double zero. Therefore there exists $t^- \in {\bI}_0$ such that 
\begin{eqnarray} 
\label{first-t}
\Theta_{\mu, \nu^-}(t^-)=\sup_{t\in {\bI}_0}{\Theta}_{\mu,\nu^-}(t) = 0 
=  {\bf P}_0(\mu,\nu^-) + \sup_{t\in {\bI}_0}{\hat {\Theta}}_{\mu,\nu^-}(t) .
\end{eqnarray}
At the last tangency, ($\nu=\nu^+$) the graph of the splitting function 
is located above the $t$ axis  and $\Theta_{\mu,\nu^+}$
admits  a double zero. 
Therefore there exists $t^+ \in {\bI}_0$ such that
\begin{eqnarray} 
\label{last-t}
\Theta_{\mu, \nu^+}(t^+)=\inf_{t\in {\bI}_0}{\Theta}_{\mu,\nu^+}(t) = 0 
=  {\bf P}_0(\mu,\nu^+) + \inf_{t\in {\bI}_0}{\hat {\Theta}}_{\mu,\nu^+}(t) .
\end{eqnarray}
If we neglect ${\cal O}_2$ in (\ref{real-exp}),
(\ref{first-t}) and (\ref{last-t})
are equivalent to
\begin{eqnarray}
\label{mm-inter}
\left\{
\begin{array}{ccccc}
{\bf P}_0(\mu,\nu^+) 
 & - &  
\ \  2|{\bf P}_{-1}(\mu, \nu^+)| & = & 0 ,
\\  & & & & \\
{\bf P}_0(\mu,\nu^-) 
 &  + & 
2|{\bf P}_{-1}(\mu, \nu^-)|& = & 0.
\end{array}
\right.
\end{eqnarray}
In this way the problem of finding 
the first and 
the last
tangencies,  
is replaced by scalar equations in one variable each.
Therefore, instead of looking 
for
intersections between $W^u_{\mathrm {loc}}$ and $W^u$
and their tangencies,
 we save a lot of time by simply solving  a scalar equation.
Observe that for $\nu$ near ${\bar \nu}$, for all $t \in {\bI}_0$ 
 we have
\begin{eqnarray}
\label{p-o-p-1}
{\bf P}_0(\mu, \nu)  =  
{\cal O}(|{\bf P}_1(\mu, \nu)), 
\qquad \ \sup_{t \in {\bI}_0}|{\Theta}_{\mu,\nu}(t)| =  {\cal O}(|{\bf P}_1(\mu, \nu)).
\end{eqnarray}
From (\ref{w-w}) we need only 4 points per-period 
to evaluate ${\bf P}_0$ and ${\bf P}_{\pm 1}$. 
	Concretely we write
\noindent
\begin{eqnarray}\label{harmonics}
 \left\{
 \begin{array}{rcl} 
  {{\gmp}}_0(\mu,\nu) \approx {\ma}_0(\mu,\nu) & = &
 \frac{1}{2}({\tilde \Theta}_{\mu,\nu}(0)+{\tilde \Theta}_{\mu,\nu}(1/2))  \\ 
& \ & \\
  {{\gmp}}_{-1}(\mu,\nu)   \approx {\ma}_{-1} (\mu,\nu)
 & = & \frac{1}{4}({\tilde \Theta}_{\mu,\nu}(0)-{\tilde \Theta}_{\mu,\nu}(1/2)
\\ & \ & \\
& + & i({\tilde \Theta}_{\mu,\nu}(1/4) -{\tilde \Theta}_{\mu,\nu}(-1/4))) \\
& \ & \\
 {{\gmp}}_{1}(\mu,\nu)  \approx {\ma}_1(\mu,\nu) 
  & = & \frac{1}{4}({\tilde \Theta}_{\mu,\nu}(0)-{\tilde \Theta}_{\mu,\nu}(1/2) 
\\ & \ & \\
& -& i({\tilde \Theta}_{\mu,\nu}(1/4) -{\tilde \Theta}_{\mu,\nu}(-1/4))). 
 \end{array}\right. \end{eqnarray}
From (\ref{w-w}) and (\ref{formula}), the approximation here means
\begin{eqnarray}\label{ror}
\max \{ |{\bf R}_0(\mu,\nu) - {\bf P}_0(\mu,\nu)| , 
|{\bf R}_{\pm 1}(\mu,\nu) - {\bf P}_{\pm 1 }(\mu,\nu)|\}
={\cal O}(\sup_{t \in {\bI_0}}|\Theta_{\mu,\nu}(t)|^2). 
\end{eqnarray}
Moreover, 
with (\ref{p-o-p-1})
we have
\begin{eqnarray}
\label{t-w-have}
|{{\ma}}_{\pm 1}(\mu,\nu )-{{\gmp}}_{\pm 1}(\mu,\nu ) | & = & 
{\cal O}(|{{\ma}}_{\pm 1}(\mu,\nu )|^2).
\end{eqnarray}
We then solve
\begin{equation}
\label{avec-tilde}
\left\{
\begin{array}{ccccc}
{\bf R}_0(\mu,{\tilde \nu}^+) 
  & - &  
\ \  2|{\bf R}_{-1}(\mu, {\tilde \nu}^+)| & = & 0 ,
\\  & & & & \\
{\bf R}_0(\mu,{\tilde \nu}^-) 
 & + &  2|{\bf R}_{-1}(\mu, {\tilde \nu}^-)|& = & 0.
\end{array}
\right.
\end{equation}
From (\ref{p-o-p-1}), (\ref{ror}),  (\ref{t-w-have}) and (\ref{avec-tilde}), we have
\begin{eqnarray}
\label{almost}
\Theta_{\mu,{\tilde \nu}^-}(t^-) =  {\cal O}({\bf R}_{-1}^2(\mu, {\tilde \nu}^- ) ), \qquad 
\Theta_{\mu,{\tilde \nu}^+}(t^+) =  {\cal O}({\bf R}_{-1}^2(\mu, {\tilde \nu}^+ ) ).  
\end{eqnarray}
By the Mean Value Theorem, we have
\begin{eqnarray}
\label{e-w-1}
|\nu^+ -{\tilde \nu}^+ | =  
{\cal O}\biggl(\displaystyle\frac{{\bf R}_{-1}^2(\mu, {\bar \nu})}{\partial\Theta_{\mu, {\nu} }
/\partial \nu|_{\nu={\bar \nu}}
}\biggr), \qquad  
|\nu^- -{\tilde \nu}^-|  =  
{\cal O}\biggl(\displaystyle\frac{{\bf R}_{-1}^2(\mu, {\bar \nu})}{\partial\Theta_{\mu, {\nu} }
/\partial \nu|_{\nu={\bar \nu}}
}\biggr). 
\end{eqnarray}
%


\vskip10pt

This approach 
gives a good estimation
of the locus of the homoclinic 
zone and therefore of the corresponding width, but requires
the computation of 
both $\nu^+$ and $\nu^-$ with  a very high precision.
To be more precise, 
assume we want to compute the width of the homoclinic
zone for a given 
value of the main parameter with $N$ correct digits, 
while the width of the zone (roughly estimated with formula (6)) satisfies
\begin{eqnarray}
\label{z-10}
10^{N_z+1} \leq \nu^+ -\nu^- < 10^{N_z},
\end{eqnarray}
where $N_z>>1$. Thus we need to compute both $\nu^+$ and $\nu^-$
with $N_z+N$ correct digits. We observe (numerically) that
\begin{eqnarray}
\label{v-anti-1}
{\tilde \nu}^+ - {\tilde \nu}^- & = & {\cal O}\biggl(
\displaystyle\frac{|{\bf R}_{-1}|(\mu,{\bar \nu})}{\partial \Theta_{\mu,{\nu}}(t_0)|_{\nu={\bar \nu}}}
\biggr),
\end{eqnarray}
also compare with
(\ref{width-exp}) below.
Therefore with (\ref{e-w-1}) and (\ref{v-anti-1}),
${\tilde \nu}^+ -{\tilde \nu}^-$ gives 
an estimation of  
the width with a relative error of the same order. 
In particular, this  means that we cannot choose $N$  bigger than $N_z$.
With this method, thanks to (\ref{t-w-have}), 
the estimations of  ${\bf P}_0(\mu,\nu)$ and of
${\bf P}_{-1}(\mu,\nu)$  are obtained
with a relative error of the same order 
as $|{\bf P}_{-1}(\mu, {\bar \nu})|$.
This 
requires the  computation of the splitting determinant 
with the same relative precision.
When the main parameter tends to $0$, 
since the eigenvalues $\lambda_1$ and $\lambda_2$ tend to $1$, 
the number of iterations  (i.e., $m_0$) 
 and the number of terms in 
(\ref{ust}), (i.e., $N_{\mathrm {max}}$)  required to compute  
the unstable manifold  need to be chosen bigger and bigger. 
Moreover, in order to guarantee (\ref{t-w-have}), 
we need to have 
${\bf P}_0(\mu,\nu) ={\cal O}({\bf P}_{-1}(\mu,\nu))$, 
i.e., (\ref{p-o-p-1}), which requires that
the local stable and the 
unstable manifold are close to one another and more precisely
\begin{eqnarray}\label{ggk}
\|{\bf \Gamma}_u(t)-{\bf \Gamma}_u(t))\| = {\cal O}(K(\mu,\gamma-2)).
\end{eqnarray}
As a conculsion, when the main parameter 
tends to $0$, this approach becomes more and more delicate.

In what follows, we propose another approach which does not 
require the computation of ${\bf P}_0(\mu,\nu)$, still requires a first 
value of $\nu={\bar \nu}$ such that
(\ref{p-o-p-1})
and gives an estimation of 
the width with the same precision.

%
%
%

\subsection{`Real' approach}

From (\ref{first-t}) and (\ref{last-t}) we have
\begin{eqnarray}
\label{tog}
 {\bf P}_0(\mu,\nu^+) - {\bf P}_0(\mu,\nu^-) & = 
& -\inf_{t\in {\bI}_0}{\hat \Theta}_{\mu,\nu^+}(t) + 
\sup_{t\in {\bI}_0}{\hat \Theta}_{\mu,\nu^-}(t).
\end{eqnarray}
Furthermore, 
from the Mean Value Theorem, 
there 
exists 
$\nu^-\leq \nu_2 \leq \nu ^+$ 
such that 
\begin{eqnarray}
\label{mean-v-t}
 {\bf P}_0(\mu,\nu^+) - 
{\bf P}_0(\mu,\nu^-) & = & 
\displaystyle\frac{\partial {\bf P}_0}{\partial \nu}  |_{\nu =\nu_2}
\cdot (\nu^+-\nu^-).
\end{eqnarray}
Thus we get
\begin{eqnarray}
\label{width-exp}
\nu^+-\nu^- & = 
& \displaystyle\frac{\sup_{t\in {\bI}_0}{\hat \Theta}_{\mu,\nu^-}(t)-
\inf_{t\in {\bI}_0}{\hat \Theta}_{\mu,\nu^+}(t)}
{{\partial {\bf P}_0}/{\partial \nu}  |_{\nu =\nu_2}}.
\end{eqnarray}

\noindent
We observe (numerically)  
that
 ${\hat \Theta}_{\mu,\nu}$ 
does not change much with respect to $\nu$.
More precisely
for all $\nu^-\leq \nu_3 \leq \nu^+,  \nu^-\leq \nu_4\leq \nu^+$ 
and for all $t \in {\bI}_0$,
\begin{eqnarray}
\label{p11}
\frac{|{\hat \Theta}_{\mu,\nu_4}(t) -
{\hat \Theta}_{\mu,\nu_3}(t)|}{\nu_4-\nu_3} 
& = & 
 {\cal O}(|{\bf P}_{-1}|(\mu, {\bar \nu})).
\end{eqnarray}
Thus, with (\ref{real-exp}) and (\ref{p-o-p-1}) we have
\begin{eqnarray}
\label{wr}
\sup_{t\in {\bI}_0}{\Theta}_{\mu,\nu^-}(t)- 
\inf_{t\in {\bI}_0}{\Theta}_{\mu,\nu^+}(t) & = &   
4|{\bf P}_{-1}|(\mu, {\bar \nu}) +{\cal O}({\bf P}_{-1}^2(\mu, {\bar \nu})). 
 \end{eqnarray}

\noindent
Furthermore with  (\ref{hat-theta}) we have
\begin{eqnarray}
\label{p-0-theta}
\displaystyle\frac{{\partial {\bf P}_0}}{{\partial \nu}}  |_{\nu =\nu_2}(t) 
& = &
\displaystyle\frac{\partial \Theta_{\mu,\nu}}{\partial \nu}  |_{\nu =\nu_2}(t)   
- 
\displaystyle\frac{{\partial {\hat \Theta}_{\mu,\nu}}}{{\partial \nu}}  |_{\nu =\nu_2}(t) \\
\nonumber 
& = & \displaystyle\frac{\partial \Theta_{\mu,\nu}}{\partial \nu}  |_{\nu =\nu_2}(t) 
+{\cal O}({\bf P}_{-1}(\mu,{\bar \nu})).
 \end{eqnarray}
With (\ref{def-tet}) and (\ref{formula}) we have  
\begin{eqnarray}
\label{t-tilde}
\displaystyle\frac{{\partial { \Theta}_{\mu,\nu}}}{{\partial \nu}} 
 |_{\nu =\nu_2}(t)  =  
\displaystyle\frac{{\partial {\tilde \Theta}_{\mu,\nu}}}{{\partial \nu}}  |_{\nu =\nu_2}(t) +
{\cal O}({\tilde \Theta}_{\mu,{\bar \nu}}(t)).
\end{eqnarray}
We then  write
\begin{eqnarray}
\label{f0}
{\tilde \Theta}_{\mu,\nu_4}(t) - {\tilde \Theta}_{\mu,\nu_3}(t) &  = & 
\displaystyle\frac{\partial  {\tilde \Theta}_{\mu,\nu}}{\partial \nu}(t)  |_{\nu =\nu_3} \cdot
(\nu_4-\nu_3) +{\cal O}(({\nu_4-\nu_3})^2),
\end{eqnarray}
and therefore

\begin{eqnarray}
\label{f01}
\displaystyle\frac{{\tilde \Theta}_{\mu,\nu_4}(t) - {\tilde \Theta}_{\mu,\nu_3}(t)}{\nu_4-\nu_3}
& = & \displaystyle\frac{\partial  {\tilde \Theta}_{\mu,\nu}}{\partial \nu}(t)  |_{\nu =\nu_2} 
+{\cal O}(({\nu^+-\nu^-})).
\end{eqnarray}
We observe (numerically) 
 that the left hand side of (\ref{f01}) stays away from $0$ as the main
parameter tends to $0$, more precisely there exists $v_0>0$ such that 
for all $\mu>0$, $\nu_3$, $\nu_4$ near ${\bar \nu}$ and for all $ t \in {\bI}_0$,
\begin{eqnarray}
\label{vitesse}
|\displaystyle\frac{{\tilde \Theta}_{\mu,\nu_4}(t) - 
{\tilde \Theta}_{\mu,\nu_3}(t)}{\nu_4-\nu_3}|
& > & v_0. 
\end{eqnarray}
With (\ref{p-o-p-1}), (\ref{p-0-theta}) and (\ref{f01})
and by choosing $\nu_3$ and $\nu_4$ sufficiently close to 
one another,  we have
\begin{eqnarray}
\label{f2}
\displaystyle\frac{{\partial {\bf P}_0}}{\partial \nu}  |_{\nu =\nu_2}  =   
\displaystyle\frac{{\tilde \Theta}_{\mu,\nu_4}(t) - {\tilde \Theta}_{\mu,\nu_3}(t)}{\nu_4-\nu_3} 
+
{\cal O}(|{\bf P}_{-1}|(\mu, {\bar \nu}))
+ {\cal O}(\nu^+-\nu^-).
\end{eqnarray}
Therefore, with (\ref{width-exp}), (\ref{wr}), (\ref{vitesse}) 
and (\ref{f2}) we can write 
\begin{eqnarray}
\label{zmu}
\nu^+-\nu^-  & =&  {\cal Z}(\mu) + {\cal O}({\cal Z}^2(\mu))
\end{eqnarray}
where
\begin{eqnarray}
\label{w-e-app}
{\cal Z}(\mu)& = & 
\displaystyle\frac{4|{\bf P}_{-1}|(\mu,{\bar \nu})(\nu_3-
\nu_4)}{{\tilde \Theta}_{\mu,\nu_3}(t_0) -{\tilde \Theta}_{\mu,\nu_4} (t_0)},
\end{eqnarray}
where $t_0 \in {\bI}_0$.




\noindent
Thanks to (\ref{harmonics}) 
we obtain the following estimation for the width of 
the homoclinic zone
\begin{eqnarray}
\label{z-real}
\nu^+-\nu^- \approx {\cal Z}_r(\mu) =  
\displaystyle\frac{4|{\bf R}_{-1}(\mu,{\bar \nu})|(\nu_3-
\nu_4)}{{\tilde \Theta}_{\mu,\nu_3}(t) -{\tilde \Theta}_{\mu,\nu_4} (t)}.
\end{eqnarray}
With (\ref{t-w-have}),  (\ref{zmu})  and (\ref{w-e-app})
it follows that
\begin{eqnarray}
\label{r-v-sone}
(\nu^+-\nu^-)- {\cal Z}_r(\mu)  & = &{\cal O}({\cal Z}^2_r(\mu)).
\end{eqnarray}

This 'real' approach gives a good estimation of the width of the homoclinic
zone with the same precision 
as before in (\ref{e-w-1}). Moreover, it 
requires only 
the computation of ${\bf P}_{-1}(\mu, {\bar \nu})$
and that of ${\tilde \Theta}_{\mu, \nu}(t)$ for two different
values
of $\nu$. However, we still need to find a value of $\nu={\bar \nu}$
such that (\ref{p-o-p-1}) holds.


\vskip10pt
In what follows we  
present another  way to compute the width: in the new approach, 
${\bf \Gamma}_u$ does not need to return near ${\bf \Gamma}_s$ as close as in
(\ref{ggk}). In this way, we will be able to compute   
the splitting determinant  
with less
precision.
This alternative approach 
consists of looking at the splitting function for complex
value of $t$.

\subsection{'Complex' approach}
Now we present another way to compute the first harmonic, with less precision
than in the 'real' approach case, but with less effort. 
Recall that  
formulae (5) and (6) already give the following estimate 
\begin{eqnarray}
\label{z-m-o}
\nu^+ -\nu^-  
= {\cal O} ( K(\mu, \gamma -2)). 
\end{eqnarray}
Moreover,  with (\ref{vitesse}) 
and (\ref{w-e-app}), (\ref{z-m-o})  gives us a rough estimate 
of $|{\bf P}_{\pm 1}|$,
i.e., 
we have 
$|{\bf P}|_{ \pm 1}(\mu, \nu) ={\cal O} (K(\mu, \gamma-2))$. 

\vskip10pt
Take $0<<\delta<\varrho$ and 
$\Delta_0 = K(\mu, \gamma -2)e^{2\pi\delta}$
such that
$K(\mu, \gamma -2)<<\Delta_0 $.
Assume that we have found a value of $\nu=\nu_0$ such that
\begin{eqnarray}
\label{1-ap}
K(\mu, \gamma-2)<< \sup_{t\in {\bI}_0}\Theta_{\mu,\nu_0 } (t) 
\leq  \Delta_0.
\end{eqnarray}
Observe that looking for such a value of $\nu= \nu_0$ requires 
less effort than searching 
for ${\bar \nu}$ where 
$\sup_{t\in {\bI}_0}\Theta_{\mu, {\bar \nu}} ={\cal O}(K(\mu, \gamma -2))$.
In particular, we only need to compute the splitting function
with a relative error of order $\sup_{t\in {\bI}_0}\Theta_{\mu,\nu_0 } (t)$. 
With (\ref{w-w}), there exists $s \in {\bI}_0$ such that
\begin{eqnarray*}
\sup_{t\in {\bI}_0}\Theta_{\mu, \nu_0}(t)
 = 
{\bf P}_0(\mu,\nu_0) +  
{\bf P}_{-1}(\mu,\nu_0)e^{-2i\pi s}  + {\bf P}_{1}(\mu,\nu_0)
 e^{2i\pi s} +{\cal O}_2(s).
\end{eqnarray*}
Since $|{\bf P}_{\pm 1}|(\mu,\nu_0) ={\cal O}(K(\mu, \gamma -2))$, 
we then conclude that
${\bf P}_0(\mu,\nu_0) ={\cal O}(\Delta_0)$. Since 
$\Delta_0 >> K(\mu,\gamma -2)$, 
we are not able to compute precisely 
the first 
harmonic ${\bf P}_{-1}(\mu, \nu_0)$, with the real approach. 
However, instead of considering $t \in {\bI}_0$ 
as real, we now 
consider $t$ in the complex interval $ [\delta i , \delta i +1]$. 
Recall that the Fourier 
coefficients of ${\Theta}_{\mu, \nu_0}(t)$ are
\begin{eqnarray*}
{\bf P}_0(\mu,\nu_0)  = 
 \int_{0}^1 { \Theta}_{\mu,\nu_0}(t)dt,  \ {\bf P}_{-1}(\mu,\nu_0)
 =  
\int_{0}^1 e^{2\pi i t} { \Theta}_{\mu,\nu_0}(t)dt.
\end{eqnarray*}
Since ${ \Theta}_{\mu,\nu_0}$ is periodic and analytic in ${\mathrm {B}}$
defined in (31), 
we have
\begin{eqnarray}
\label{f-c-comp-ex}
{\bf P}_{-1}(\mu,\nu_0) & = &   
\int_{i\delta}^{i\delta + 1} e^{2\pi i t} {\Theta}_{\mu,\nu_0}(t)dt. 
\end{eqnarray}
With (\ref{w-w}) we have
\begin{eqnarray}
\label{nu-ici}
e^{2\pi i t} { \Theta}_{\mu,\nu_0}(t)  =    
{\bf P}_{-1}(\mu,\nu_0) +e^{2\pi i t} {\bf P}_{0}(\mu,\nu_0)
+  e^{4\pi i t} {\bf P}_{1}(\mu,\nu_0)
 + e^{2\pi i t}{\cal O}_2(t),
\end{eqnarray}
where 
\begin{eqnarray}
\label{ineg-6}
{\cal O}_2(t) & = & 
{\cal O}(\sup_{t^\prime \in {\bI}_0}|{ \Theta}_{\mu,\nu_0}^2(i\delta+ t^\prime)|).
\end{eqnarray}
With (\ref{h-k}) we have 
\begin{eqnarray*}
{\bf P}_{\pm 1} (\mu,{\nu_0}) & = &  {\cal O}(e^{-2\pi\varrho}).
\end{eqnarray*}
Therefore, since ${\bf P}_0(\mu, \nu_0) ={\cal O}(\Delta_0)$, with (\ref{w-w}) and (\ref{ineg-6}), 
we have
\begin{eqnarray}
\label{cause}
\sup_{t^\prime \in {\bI}_0}|{ \Theta}_{\mu,\nu_0}|(i\delta+ t^\prime)  =   
{\cal O}(e^{2\pi \delta-2\pi \varrho})=
{\cal O}(|{\bf P}_{-1}|e^{2\pi \delta}),
\end{eqnarray}
and further we have
\begin{eqnarray}
 \label{ineg-1} 
 {\bf P}_{1}(\mu,\nu_0) e^{4 i \pi t }  & =  & {\cal O}(e^{-2\pi (\rho +2\delta)}), \\ 
\nonumber
\ |e^{2\pi i t}{\cal O}_2(t)| & = &  
{\cal O}(|e^{+2\pi it} |{\bf P}_{-1}|^2(\mu, \nu_0) e^{4\pi\delta})|) = 
{\cal O}(e^{-2\pi(2\rho -\delta )}).
\end{eqnarray}
We distinguish two cases

\par\vskip10pt 
{\bf Case 1}: $\delta > \displaystyle\frac{\rho}{3}$. 
In this case, 
$2\varrho -\delta< \varrho +2\delta$ and 
from (\ref{ineg-1}) we have 
\begin{eqnarray*}
|{\bf P}_{1}(\mu,\nu_0) e^{4 i \pi t }| &  <<  & |e^{2\pi i t}{\cal O}_2(t)|.
\end{eqnarray*}
Using (\ref{formula}) 
we write
\begin{eqnarray}\label{centre2}
e^{2\pi i t} { \Theta}_{\mu,\nu_0}(t) =e^{2\pi i t} {\tilde \Theta}_{\mu,\nu_0}(t) 
+ e^{2\pi i t} {\tilde h}_{\mu,\nu_0}(t) =
{\cal A}(t)+ {\cal E}_1(t)
\end{eqnarray}
where 
$
{\cal A}(t) =  {\bf P}_{-1}(\mu,\nu_0) +e^{2\pi i t } {\bf P}_0(\mu,\nu_0)
$
and  with (\ref{ineg-1}),
$$
{\cal E}_1(t) = {\cal O}(e^{-2\pi(2\rho -\delta)}).
$$
Observe that for $t \in [i\delta , i\delta +1]$
\begin{eqnarray*}
|e^{2 i \pi t }{\tilde h}_{\mu,\nu_0}(t)|
=  {\cal O}(e^{-2\pi(2\rho-\delta)}).
\end{eqnarray*}
In this case
\begin{eqnarray}\label{centre1}
\int_{i\delta}^{i\delta +1} {\cal A}(t)dt =\displaystyle\frac{1}{2}\biggl( {\cal A}(i\delta ) 
+{\cal A}(i\delta + 1/2)\biggr).  
\end{eqnarray}
But with, (\ref{centre2}) we have 
\begin{eqnarray}\label{centre3}
\int_{i\delta}^{i\delta +1} e^{2\pi i t} {\tilde \Theta}_{\mu,\nu_0}(t)dt= 
\int_{i\delta}^{i\delta +1} {\cal A}(t)dt +  {\cal O}(e^{-2\pi(2\rho-\delta)}).
\end{eqnarray}
Finally from
(\ref{f-c-comp-ex}), (\ref{centre2}), (\ref{centre1}) and (\ref{centre3})
we get
\begin{eqnarray}\nonumber
{\bf P}_{-1}(\mu,\nu_0)  & = &    {\bf C}_{-1}(\mu,\nu_0) + {\tilde r}_1, \\ \nonumber
{\mathrm {where  }} \ {\bf C}_{-1}(\mu,\nu_0) & = &  
\displaystyle\frac{1}{2}e^{-2\pi \delta} 
\biggl( {\tilde \Theta}_{\mu,\nu_0}(i\delta) - 
{\tilde \Theta}_{\mu,\nu_0}(i\delta +1/2)\biggr) 
\\ \label{ca-11}
|{\tilde r}_1| &= & {\cal O}(e^{-2\pi(2\rho -\delta)}).
\end{eqnarray}

\par\vskip10pt
{\bf Case 2}: $\delta \leq \displaystyle\frac{\rho}{3}$.
In this case, from (\ref{ineg-1}) we have
\begin{eqnarray*}
|{\bf P}_{1}(\mu,\nu_0) e^{4 i \pi t }| &  \geq   & |e^{2\pi i t}{\cal O}_2(t)|,
\end{eqnarray*}
therefore we cannot neglect 
the term ${\bf P}_{1} e^{4\pi i t }$ from the integration in 
(\ref{f-c-comp-ex}). Thus we write

\begin{eqnarray*}
e^{2i\pi t} {\tilde \Theta}_{\mu,\nu_0}(t) & = &  
{\bf P}_{-1}(\mu,\nu_0) +e^{2\pi i t } {\bf P}_0(\mu,\nu_0) +{\bf P}_{1}(\mu,\nu_0) e^{4\pi i t } 
+ {\cal E}(t) \\ \nonumber
 & = & 
{\tilde {\cal A}}(t)+ {\cal E}(t)
\end{eqnarray*}
where  with (\ref{ineg-1})
$$
{\tilde {\cal A}}(t)=  {\bf P}_{-1}(\mu,\nu_0) +e^{2\pi i t } {\bf P}_0(\mu,\nu_0) +
{\bf P}_{1} e^{4\pi i t }, 
\ {\cal E}(t) = {\cal O}(e^{-2\pi(2\rho -\delta)}).
$$
In this case
$$
\int_{i\delta}^{i\delta +1} {\tilde {\cal A}}(t)dt= \displaystyle\frac{1}{4}\biggl(
{\tilde {\cal A}}(i\delta) + {\tilde {\cal A}}(i\delta +1/2 ) + {\tilde {\cal A}}(i\delta +1/4 ) + {\tilde {\cal A}}( i\delta +3/4)
\biggr)
$$
and we get
\begin{eqnarray}\nonumber
{\bf P}_{-1}(\mu,\nu_0)  & = &    {\bf C}_{-1}(\mu,\nu_0) + {\tilde r}_2, \ {\mathrm {where}} \\ \nonumber
\ {\bf C}_{-1}(\mu,\nu_0) & = &  
e^{-2\pi\delta} 
\frac{1}{4}({\hat \Theta}_{\mu,\nu_0}(0)-{\hat \Theta}_{\mu,\nu_0}(1/2)
\\ \nonumber
 & - &  
i({\hat \Theta}_{\mu,\nu_0} (1/4) -\Theta_{\mu,\nu_0}(-1/4)))
\\ \label{ca-22}
|{\tilde r}_2| &= & {\cal O}(e^{-2\pi(2\rho -\delta)})=
{\cal O}({\bf C}^2_{-1}(\mu, \nu_0)e^{2\pi\delta}).
\end{eqnarray}

\noindent
When $\delta > \varrho/3$, the estimation given in (\ref{ca-11}) requires the
computation of the splitting 
determinant 
${\tilde \Theta}_{\mu,\nu_0}(t)$ 
at  two different values of $t$ only. However, 
when 
 $\delta \leq  \varrho/3$, (\ref{ca-22}) requires four different values of $t$.
 The computation
in the first case is faster, but since $\delta$ is bigger, we loose 
some precision.

\vskip10pt
From the 'complex' approach, 
the width of homoclinic zone is
approximated by
\begin{eqnarray} 
\label{zone-23}
\nu^+(\mu) -\nu^-(\mu)  
& \approx & 
{\cal Z}_c(\mu) 
\\  \nonumber
{\mathrm {where}} & \ &   
{\cal Z}_c(\mu)   
 =  
4\frac{|{{\hmc}}_{-1}(\mu,\nu_0)|(\nu^\prime_3-\nu^\prime_4)}
{{\tilde \Theta}_{\mu,\nu_3^\prime}(t) - {\tilde \Theta}_{\mu,\nu_4^\prime}(t)},
\end{eqnarray}
where $\nu_3^\prime$  and $\nu_4^\prime$ are chosen near $\nu_0$.
Since ${\bf C}_{-1}(\mu,\nu_0)={\cal O}(e^{-2\pi \varrho})$ 
with (\ref{ca-11}) or (\ref{ca-22}) we have
\begin{eqnarray}\label{star5} 
|{\cal Z}_c(\mu) -(\nu^+(\mu) -\nu^-(\mu))|=
{\cal O}(\hmc_{-1}^2(\mu,\nu^+)e^{2\pi\delta}).
\end{eqnarray}


\vskip10pt
\subsection{'Real' versus 'Complex'}
The real approach provides a good estimation of the width of the homoclinic 
zone. More precisely, formula 
(\ref{z-real}) 
gives an estimation 
of the width with a relative error of the same order, see (\ref{r-v-sone}).
However, this approach requires to compute the 
splitting determinant 
with  the same relative error. This task becomes more and more delicate 
as 
the main parameter approaches $0$.
The complex approach requires less  
precision for the computation of  the 
splitting determinant  (and therefore can be computed much faster)
as  $\delta$ is chosen larger. However, the 
estimation of the width is obtained with less precision. 	

In the case of the Bogdanov map, we use similar notations: 
$b$ is the slave parameter and $a$ is the main parameter.
The first harmonic computed with   (\ref{harmonics}) is denoted by 
${\bf R}_{-1}(a, {\bar b})$, where ${\bar b}$ is the analogue of ${\bar \nu}$
in the case of the Quadratic map. 
Simlarly, ${\bf C}_{-1}(a, {b}_0)$ stands for the 
first harmonic computed with (\ref{ca-11}) or (\ref{ca-22})
where ${b}_0$ is the analogue of ${\nu}_0$
in the case of the Quadratic map.
For illustration,  
we compute the first 
harmonic and the width of the 
homoclinic zone using both approaches  
for the Bogdanov map (${\tilde \gamma}=3$), see Figure 2. 
We easily verify that
$$
\log_{10}\biggl( {\hmc}_{-1}(a,b_0)-{\ma}_{-1}(a,{\bar b})\biggr) 
\approx 2\log_{10}(|{\ma }|_{-1}(a, {\bar b}))+ 
\log_{10}(e^{2\pi\delta}),
$$
which follows from (\ref{r-v-sone}) and (\ref{star5}).
Furthermore, 
we also verify that 
\begin{eqnarray*}
\log_{10}|{\cal Z}_r(a)-{\cal Z}_c(a)| & \approx  & 
 \log_{10}\biggl(
\frac{|{\hmc}_{-1}(a,b_0)|e^{2\pi\delta}(b^\prime_3-b^\prime_4)}
{{{\tilde \Theta}_{a,b^\prime_3}(t_0) - {\tilde \Theta}_{a,b^\prime_4}(t_0)}}
\biggr) \\ 
& \approx &  
\log_{10}(b^+-b^-) + \log_{10} |{\hmc}_{-1}(a,b_0)e^{2\pi\delta}|,
\end{eqnarray*}
where $b^\prime_3 $, $b^\prime_4 $ are the analogues of 
$\nu^\prime_3 $, $\nu^\prime_4 $ respectively.
\par\vskip10pt\noindent
{\bf Example}: 
We consider the Bogdanov map when $a \approx 7*10^{-5}$.
Using the real approach, we have 
$\log_{10}(b^+-b^-) \approx -1000$, see Figure 2.
With this approach, we  
compute ${\tilde \Theta}_{a, {\bar b}}$ with a relative error
of order $10^{-1000}$, which is already  a quite delicate task.
However, from the complex approach, we can 
(for instance) choose $\delta$ in such a way that
$e^{2\pi\delta} \approx 10^{700}$, see Figure 2.
This way, for values of $t \in [i\delta, i\delta+1]$, we have
$$
\log_{10}({\tilde \Theta}_{a,b_0}(t)) 
\approx \log({\bf C}_1(a, b_0) e^{2\pi\delta}) \approx
-300.$$ 

\noindent
Therefore,
computing
${\bf C}_1(a,b_0)$
with (\ref{ca-22}) requires the computation 
of the splitting determinant
with a relative error or order $10^{-300}$. 
Moreover, 
we just need to find a first value of $b=b_0$ such that
$$
\log_{10} \sup_{t\in {\bI}_0}\| {\bf \Gamma}_u(t)- {\bf \Gamma}_s(t)\| 
\approx
-300.
$$
However, instead of having a relative error 
for the width of 
order $10^{-1000}$ 
as in real approach case, we obtain an estimation of the width
 with a relative error
of order $10^{-300}$.
\begin{figure}\label{55}
\begin{center}
\includegraphics[scale=0.40]{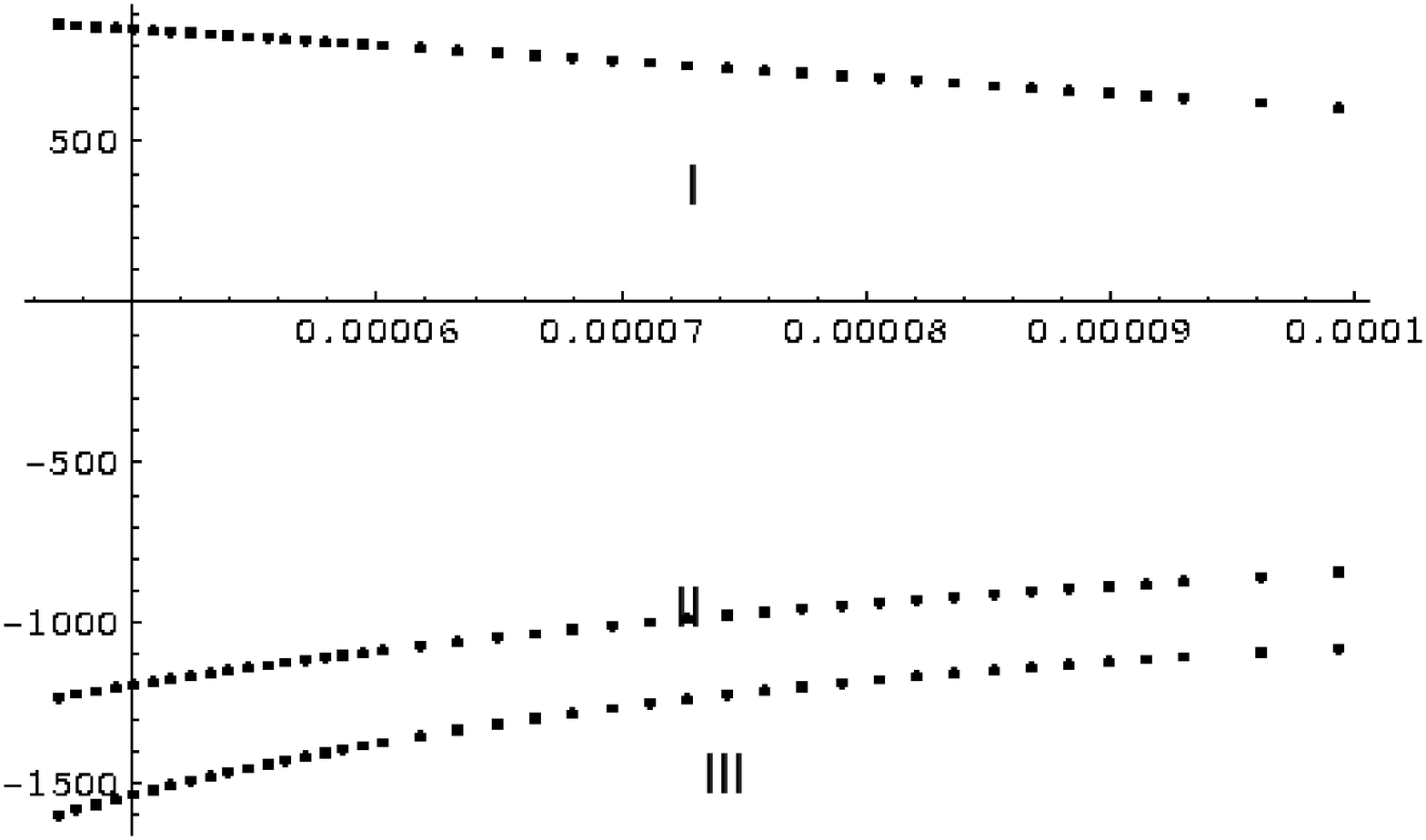}
\includegraphics[scale=0.40]{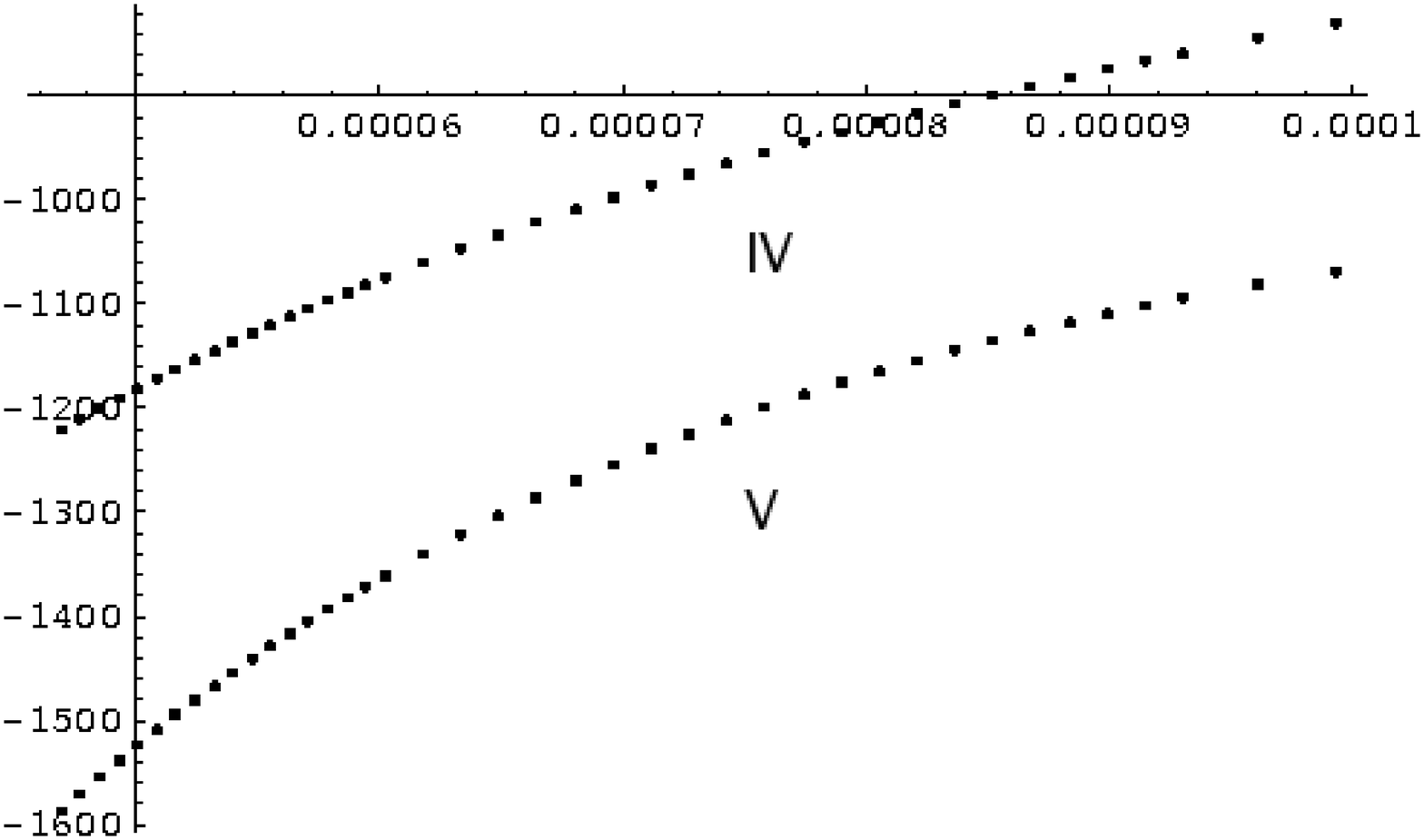}
\caption{Above: (I)-Graph of 
 $\log_{10}|e^{2\pi \delta}|$
against the parameter $a$. 
(II)-Graph of 
$\log_{10}|\hmc_1(a,b_0)|$ 
against the parameter $a$
which essentially coincides 
with the graph of 
$\log_{10}|\ma_1(a,{\bar b})|$ against $a$.  
(III)-The corresponding error i.e., 
$\log_{10}|\hmc_1(a,b_0)-\ma_1(a, {\bar  b}|$ against $a$.
Below: (IV)-Computation of the magnitude of the homoclinic zone 
with the real approach
(i.e., 
$\log_{10}({\cal Z}_r(a))$ against $a$. 
(V)-The graph of $\log_{10}|{\cal Z}_r(a)-{\cal Z}_c(a)|$ against $a$.}
\end{center}
\end{figure}

\par\vskip10pt
Now that we can compute the width of the homoclinic zone, we do so for
${\tilde n} $ 
(several hundred) values of $\mu^{1/4}$ and establish the set
\begin{eqnarray}\label{eq17}
{{\cal H}} =\{ (\mu_i^{\frac{1}{4}}, \log(\delta_i )),
\ \delta_i =\nu^+(\mu_i)-\nu^-(\mu_i),   
\ c<\mu_i< d, \ \ i=1,\ldots,{\tilde n}\}.
\end{eqnarray}
In what follows, we describe how from the ansatz  (\ref{teo-log})
we extract the corresponding coefficients.

\subsection{Extracting the coefficients}
Recall that the  ansatz 
we shall consider takes the form
(\ref{def-b}) where the $f_n$'s satisfy (\ref{dul}).
From the set ${\cal H}$ defined by (\ref{eq17}) we 
construct the following
matrices
\begin{eqnarray*}
{\bf A}= (A_{i,j})_{i=0,\ldots,{\tilde n}-1,\  j=1,\ldots,{\tilde n}}, 
\ A_{i,j} = 
f_i(\mu_j^{1/4}).
\end{eqnarray*}
In the case
of the  Bogdanov map, the set of normalised data is defined in
(\ref{bog-nom}), that is the $\mu_i^{1/4}$'s above are 
replaced by $a^{1/2}$.

\vskip10pt
Let
$$
{\bf \alpha} = (\alpha_1,\ldots,\alpha_{{\tilde n}})=
{\bf A}^{-1}\cdot {\bf w}, 
$$
${\mathrm {where}} \ 
{\bf w}= (\log \delta_1,\ldots,\log \delta_{{\tilde n}}).
$
Observe that 
\begin{eqnarray*}
\sum_{i=0}^{{\tilde n}-1}\alpha_i f_i(\mu_j^{1/4}) =
\log \delta_j, 
\ \forall j=1,\ldots,{\tilde n},
\end{eqnarray*}
that is the 
coefficients 
$\alpha_i$'s have been constructed in such a way that
the map
\begin{eqnarray}\label{d2}
{ {\phi}}^{\{{\tilde n}\}} & : & (0,\varepsilon_0)\rightarrow {\bR}, 
\ x \mapsto { \phi}^{\{{\tilde n}\}}
(x)=\sum_{i=0}^{{\tilde n}}\alpha_i f_i(x)
\end{eqnarray}
interpolates
the set of data 
${{\cal H}}$.

\par\vskip10pt
To illustrate our techniques, Table 1 indicates
the first coefficients of the interpolation (${\tilde n} \approx 100$)
in the case of the Bogdanov map (left, ${\tilde \gamma}=3$) and 
in the case of 
 the 
Quadratic map  (right, ${ \gamma}=-3$).
In the case of the H\'enon map, 
replacing the ansatz (\ref{teo-log}) by (\ref{bog-diverge}),
we obtain the coefficients indicated in Table 2.

\vskip10pt
Redoing the above interpolation for different values of $\gamma$ 
reveals that 
the first non linear terms in the 
expansion 
satisfies
\begin{eqnarray}\label{qq}
N_1({\gamma}) =-\displaystyle\biggl(\frac{6{(\gamma-2) } }{7\sqrt{2}}\biggr)^2, 
\end{eqnarray} 
in the case of the Quadratic map, 
and 
\begin{eqnarray}\label{bb}
B_1({{\tilde \gamma}}) =-\displaystyle\biggl(\frac{6{\tilde \gamma} }{7}\biggr)^2,
\end{eqnarray}
in the case of the Bogdanov map.
These equalities are  verified with a large precision. More precisely, 
we show that 
(\ref{qq}) and  (\ref{bb}) are verified up to the same number of correct 
digits as 
in (\ref{fi}) when checking 
the extrapolation to zero, see section 4.3
for more details.

\section{Validation of numerical method}
To test the validity of our result, 
 we propose three tests.
To begin with, we test the validity of the ansatz.
In what follows the experiments are presented in the cases of the Bogdanov map and the
H\'enon map, but the same test can be applied
in the case of the Quadratic map hereby confirming formula (\ref{s-r-diverge}). 

\subsection{Extrapolability}
We claim that the ansatz (\ref{teo-log}) 
is appropriate for an asymptotic expansion of the width 
if the following criterion is satisfied.

\vskip10pt
Assume a function $G: (0,\varepsilon_0) \rightarrow {\bR}$, possesses 
the following asymptotics at $0$
\begin{eqnarray*}
G(x) & \asymp & \sum_{i=0}^{\infty }\alpha_i f_i(x)
\end{eqnarray*}
where $\{ f_i(x), \ i\in {\bN} \}$
is the asymptotic sequence defined 
in (\ref{dul}).
 Define
\begin{eqnarray}\label{part}
G^{\{3k+3\}}(x) = \sum_{i=0}^{3k+3 }\alpha_i f_i(x).
\end{eqnarray}
We have 
\begin{eqnarray}\label{bef-log}
|G(x)-G^{\{3k+3\}}(x)| = x^{2k+3}\biggl(\alpha_{3k+4}+\varepsilon_1 (x)\biggr),
\end{eqnarray}
where $\varepsilon_1(x) ={\cal O}(x).$
From (\ref{bef-log}) we get
\begin{eqnarray}
\nonumber
\log |G(x)-G^{\{3k+3\}}(x)| & = & \log |\alpha_{3k+4}| +(2k+3)\log x + 
\log \biggl(1+\varepsilon(x)\biggr) \\ 
\label{after-log}
& = & \log |\alpha_{3k+4}| +(2k+3)\log x +\varepsilon_2 (x)
\end{eqnarray}
where $\varepsilon_2(x) ={\cal O}(|x|).$
This implies that the quantity
$\log|G(x)-G^{\{3k+3\}}(x)|$ 
is approximatively linear in $\log x$. 
This must be satified for
values of $x$  outside the data set used for interpolation.

\vskip10pt
Now we apply this criterion to the Bogdanov family. Recall that
$a$ is the slave parameter and $b$ is the main parameter.
Take an interval $[c^\prime, d^\prime]$ where
$c<c^\prime<d^\prime<d$ and 
consider the interpolation 
of the set ${\tilde {\cal H}}$
for values
of $a$ in $[c^\prime, d^\prime]$. 
In other words  
we consider the set 
\begin{eqnarray*}
{\tilde {\cal H}}^\prime 
& = & \{ (a^{1/2},\log \delta(a)) \in {\tilde {\cal H}} \ | \ 
c^\prime< a < d^\prime
\} 
\end{eqnarray*}
that consists of $3k+4$ different values and construct 
the corresponding set of coefficients $\{\alpha_i \}_{i=0,\ldots,3k+3}$
as described in section 3.
We plot  the set
\begin{eqnarray}\label{set-inter}
{\cal L}_{c,d}= \{(\log(a), \log|{G}^{\{{3k+3}\}}(\sqrt{a})-(b^+(a) -b^-(a))|),
\ c<a<d \} 
\end{eqnarray}
in Figure 3: ${\tilde n}=140$, 
$k=36$, $c=3.5*10^{-5}$, $d=9.4*10^{-3}$. 
The bold line shows the interval  $[c^\prime, d^\prime]$.
From (\ref{part}) and (\ref{bef-log})
we must get
$$
\log|{G}^{\{{3k+3}\}}(\sqrt{a})-(b^+(a) -b^-(a))|
\approx 37 \log a + {\mathrm {C}.} = 74 \log \sqrt{a} + {\mathrm {C}},
$$
where ${\mathrm {C}}$ is a constant.
In Figure 3, 
the set (\ref{set-inter}) 
looks like  a straight line
with a slope 
$\approx 75$, which indicates 
that the ansatz (\ref{teo-log})  satisfies the above criterion.

\begin{figure}\label{f3}
\begin{center}
\includegraphics[scale=1.35]{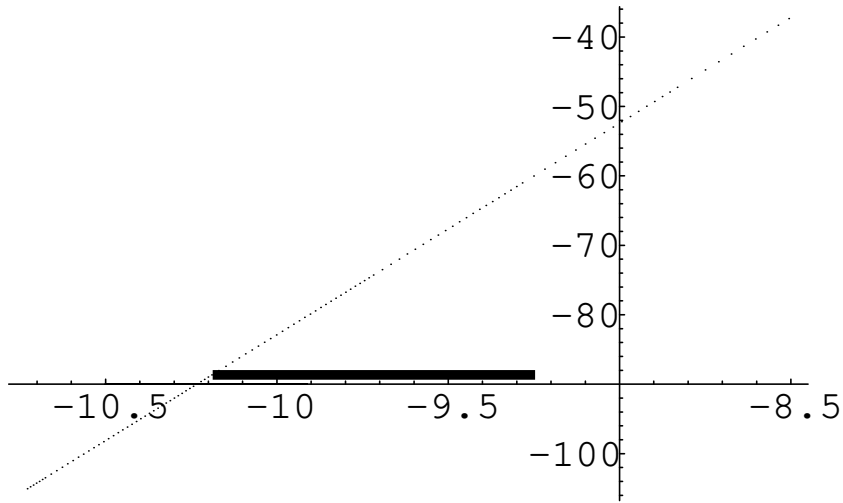}
\end{center}
\caption{
Plot of the set ${\cal L}_{c,d}$ (\ref{set-inter}) for the Bogdanov map in the case
${\tilde \gamma}=3$, $c=3.5 * 10^{-5}$, $d=9.4*10^{-3}$, ${\tilde n} =140$, $k=36$.
}
\end{figure}
\vskip10pt

In the H\'enon case, 
we interpolate the 
set of data (\ref{t-c-z}) 
with the polynomial ansatz (\ref{polyn})
and 
the normalised width
takes the form
\begin{eqnarray*}
\frac{ {\tilde b}^+({\tilde a})- {\tilde b}^- ({\tilde a}) }{ K(1-{\tilde a},0) } 
 \asymp \sum_{i=0}^{\tilde n} 
{\tilde A}_i
(1-{\tilde a})^{i/4}.
\end{eqnarray*}
We test the polynomial expansion the same way
 we test the Dulac expansion for the
Bogdanov.
More precisely, writing
$$
{\tilde G}^{\{{\tilde k}-1\}}(x) =\sum_{i=0}^{{\tilde k}-1} {\tilde A}_i x^i, \ {\tilde G}(x) 
\asymp \sum_{i=0}^{\infty} {\tilde A}_i x^i
$$
we have
\begin{eqnarray}\label{pol}
\log|{\tilde G}(x)-{\tilde G}^{\{{\tilde k}-1\}}(x)| = \log |{\tilde A}_{{\tilde k}}|+  
{\tilde k} \log x +{\cal O}(x)
\end{eqnarray}
and replacing $x$ by $({\tilde a}-1)^{1/4}$ in (\ref{pol}) leads to
\begin{eqnarray*}
\log|{\tilde G}(({\tilde a}-1)^{1/4})-{\tilde G}^{\{{\tilde k}\}}
(({\tilde a}-1)^{1/4})|  & = &
 \log |{\tilde A}_{\tilde k}|+  
\frac{{\tilde k}}{4} \log ({\tilde a}-1))
\\ \nonumber
 & + & {\cal O}(({\tilde a}-1)^{1/4}).
\end{eqnarray*}
The set 
\begin{eqnarray}\nonumber
{\tilde {\cal L}}_{c,d} 
& =&  
 \{(\log({\tilde a}-1), 
\log|{{\tilde G}}^{\{{{\tilde k}}\}}(({\tilde a}-1)^{1/4})-({\tilde b}^+({\tilde a}) -{\tilde b}^-
({\tilde a}))|), \\ \label{set-h}
& \ & 
\ c < {{\tilde a}-1} <d \}. 
\end{eqnarray}
is plotted 
(with ${\tilde k}=60$, $c=1.69 * 10^{-10}$, $d=1.125*10^{-7}$)
in Figure 4 and mimics a straight line of slope $\approx 15 =60/4$, meaning that
the polynomial ansatz satisfies the above criteria.
%



\begin{figure}[ht33333]
\includegraphics[scale=1.35]{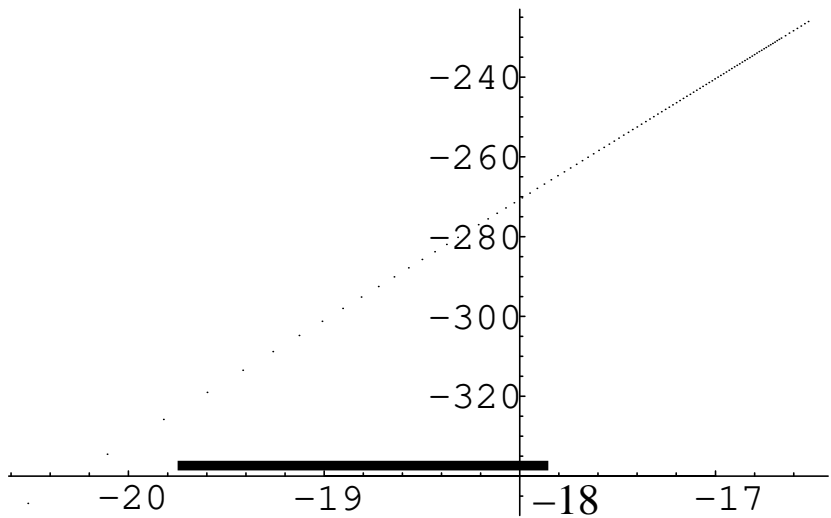}
\caption{
Plot of the set ${\tilde {\cal L}}_{c,d}$ (\ref{set-h}) for the H\'enon map, 
$c=1.69 * 10^{-10}$, $d=1.125*10^{-7}$, 
${\tilde k}=60$, ${\tilde n}= 140$.
}
\end{figure}

\vskip20pt
\noindent

The second experiment consists of checking the stability 
of our interpolation when changing (randomly) 
the data ${\tilde {\cal H}}$.

\subsection{Checking numerical stability}
In this section, our interest is with the precision  
of our data for the normalised width of the 
homoclinic zone that is required 
 in order to produce reliable results for the coefficients. 
The result of our test is presented in the 
case of the Bogdanov map, i.e., we test 
the asymptotic expansion (\ref{bog-diverge-henon}).
In order to simulate round-off errors,
we modify  the data in the $N$-th digit 
by adding a random perturbation of order $10^{-N}$ to every value of
the normalised width and recompute the coefficients of the asymptotic
expansion using the procedure described in section 3.
We repeat the experiment for several values
of $N$. Figure 5 concerns 
the coefficients $A_{11}$ in (\ref{bog-diverge-henon}): for each value of
$N$, we recompute the corresponding coefficient (denoted by 
$A_{11}^{\{N\}}$)
after adding 
a random perturbation of order $10^{-N}$.
Figure 5
clearly indicates, 
that the precision of the computation  
decrease linearly 
with respect
to $N$ and 
if the number of correct digits in the data is less than 
$N=170$, 
the corresponding 
coefficient cannot be computed correctly.
However, if $N=200$ the corresponding 
coefficient is computed with 
30 correct digits.

\begin{figure}\label{5i5}
\begin{center}
\includegraphics[scale=1.40]{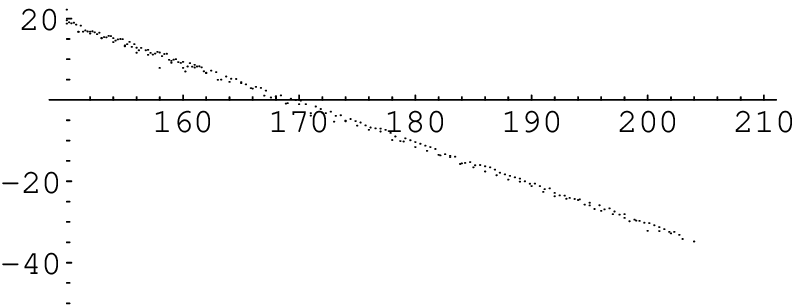}
\caption{
The relative error 
$\displaystyle\log_{10}| \frac{A_{11} -A_{11}^{\{N\}})}{A_{11}}|$
plotted against $N$ for the Bogdanov map with ${\tilde \gamma}=3$.
%
%
}
\end{center}
\end{figure}

\subsection{Extrapolation to zero}
As announced in Section 2, for each family we are able to 
define the splitting constant 
associated with the 
`unperturbed map', see \cite{GN}  for more details.
In what follows, our discussion concerns the Bogdanov family.
The splitting constant is denoted by ${\mathbf{\Theta}}({\tilde \gamma})$.
Using formula (\ref{bog-diverge-henon}), we have
\begin{eqnarray}
\label{fi}
\exp(A_0({\tilde \gamma})) & = & {\mathbf{\Theta}}({\tilde \gamma}).
\end{eqnarray}
Since we can independently compute the invariant ${\mathbf{\Theta}}({\tilde \gamma})$ 
with a very high precision, 
we can easily check the validity of our computation 
for the first term  of the asymptotic expansion.
The following table indicates, for different value of ${\tilde \gamma}$
the values of ${\mathbf{\Theta}}({\tilde \gamma})$ (left) 
computed with 20 correct digits.
For each value of ${\tilde \gamma}$, we observe that (\ref{fi})
holds and we indicate the 
 relative error represented by
 $-\log_{{\mathrm {10}}}
|({\mathbf{\Theta}}({\tilde \gamma})-\exp(A_0({\tilde \gamma})))/{\mathbf{\Theta}}({\tilde \gamma})|$
in the right column. 

\newpage
\vskip20pt
\begin{table}
\begin{tabular}{|c|c|c|}
\hline  ${\tilde \gamma}$ & ${\mathbf{\Theta}}({\tilde \gamma})$  &  $-\log_{{ \mathrm{10}}}
|({\mathbf{\Theta}}({\tilde \gamma})-\exp(A_0({\tilde \gamma})))/{\mathbf{\Theta}}({\tilde \gamma})| $ 
\\ \hline 
${\tt -2}$ & ${\tt 0.28524883190581352}\ \ \ \ \ \ \ \ $ &  ${\tt 65.23}$ \\\hline
${\tt 0}$ & ${\tt 2.47442559355325105*10^{6}}  \ \ \  \ \ \ \ \ $  & ${\tt 90.01}$\\ \hline
${\tt 3}$ & ${\tt 4.05522622851113044 *10^{26}}$ & ${\tt 62.04}$\\ \hline
${\tt 6}$ & $ {\tt 2.70980378082897208*10^{47}}$ & ${\tt 60.03}$ \\ \hline
${\tt 7}$ & ${\tt 3.09943158275750458*10^{54}}$ & ${\tt 59.6} $  \\ \hline
${\tt 9}$ & ${\tt 5.18377311752952789*10^{68}}$ & ${\tt 55,6}$  \\ \hline
\end{tabular} 
\caption{The value of  ${\mathbf{\Theta}}({\tilde \gamma})$  for different values
of ${\tilde \gamma}$. We clearly observe that the splitting constant
coincides with the first term in  (\ref{bog-diverge-henon})} 
up to the first 50 digits at least.
\end{table}
\par\vskip10pt

\vskip10pt
\noindent 
{\bf Acknowledgements}: 
This work is supported by the EPSRC grant
EP/COOO595/1.

\end{document}